\documentclass[11pt,oneside]{article}
\usepackage{amssymb,amsmath,amsthm}
\usepackage{url}
\usepackage[english]{babel}
\usepackage{pstricks,pst-node}
\renewcommand{\le}{\leqslant}
\renewcommand{\ge}{\geqslant}
\newcommand{\bad}{\mathbf{Bad}}

\renewcommand{\L}{{\rm L}}

\newcommand{\rad}{\mathrm{rad}}
\newcommand{\cent}{\mathrm{cent}}
\newcommand{\diam}{\mathrm{diam}}

\newcommand{\BBB}{\mathcal{B}}
\newcommand{\DDD}{\mathcal{D}}
\newcommand{\RR}{\mathbb{R}}

\newcommand{\ZZ}{\mathbb{Z}}
\newcommand{\QQ}{\mathbb{Q}}

\newcommand{\NN}{\mathbb{N}}

\newcommand{\III}{\mathcal{I}}

\newcommand{\CCC}{\mathcal{C}}
\newcommand{\KKK}{\mathcal{K}}
\newcommand{\LKK}{\mathcal{LK}}
\newcommand{\RRR}{\mathbf{R}}

\newcommand{\SSS}{\mathcal{S}}
\newcommand{\MMM}{\mathcal{M}}
\newcommand{\HHH}{\mathcal{H}}
\newcommand{\EEE}{\mathcal{E}}

\newcommand{\vx}{\mathbf{x}}

\newcommand{\vy}{\mathbf{y}}
\newcommand{\vr}{\mathbf{r}}

\newcommand{\vs}{\mathbf{s}}
\newcommand{\vd}{\mathbf{d}}
\newcommand{\vt}{\mathbf{t}}
\newcommand{\vp}{\mathbf{p}}
\newcommand{\vf}{\mathbf{f}}

\newtheorem*{MDP}{Mass Distribution Principle}
\newtheorem{lemma}{Lemma}
\newtheorem{theorem}{Theorem}
\newtheorem*{theorem6}{Theorem $\mathbf{6^*}$}
\newtheorem*{theoremKTV}{Theorem KTV}
\newtheorem*{theoremBPV}{Theorem BPV}
\newtheorem*{propositionBPV}{Proposition BPV}
\newtheorem*{theoremBLMV}{Theorem BLMV}
\newtheorem*{theoremEFS}{Theorem EFS}
\newtheorem*{theoremBV}{Theorem BV}
\newtheorem*{theoremB}{Theorem B}

\newtheorem{proposition}{Proposition}
\newtheorem{problem}{Problem}

\newtheorem*{corollary}{Corollary}
\newtheorem{corol}{Corollary}
\newtheorem*{definition}{Definition}

\newcommand{\cR}{{\cal R}}

\newcommand{\mad}{\mathbf{Mad}}

\begin{document}

\title{Cantor-winning sets and their applications}

\author{
 Dzmitry Badziahin, Stephen Harrap
\footnote{Research supported by EPSRC  Grant EP/L005204/1} }


\maketitle

\begin{abstract}
We introduce and develop a class of \textit{Cantor-winning} sets that share the same amenable
properties as the classical winning sets associated to Schmidt's $(\alpha,\beta)$-game: these include maximal Hausdorff
dimension, invariance under countable intersections with other
Cantor-winning sets and invariance under bi-Lipschitz
homeomorphisms. It is then demonstrated that a wide variety of badly
approximable sets appearing naturally in the theory of Diophantine
approximation fit nicely into our framework. As applications of
this phenomenon we answer several previously open questions, including some related to the Mixed
Littlewood conjecture and the $\times2, \times3$ problem.
\end{abstract}

\section{Introduction}

\subsection{Badly approximable sets}\label{sec_introbad}

The set $\bad$ of badly approximable
real numbers plays an important role in the theory of
Diophantine approximation. Recall, a real number $x$ is called \textit{badly
approximable} if there exists a constant $c(x)>0$ such that
$|x-p/q|\ge c(x)/q^2$ for all rational numbers $p/q$. It is well known that the set of all badly approximable numbers is very small in the sense that it has zero Lebesgue measure. However, a classical
 result of Jarn\'ik \cite{jarnik_1928} states that this set is in some sense as large as it can be in that it is of full
Hausdorff dimension, i.e. $\dim\bad = \dim\RR = 1$. In later works
of Davenport \cite{davenport_1964}, Pollington \& Velani \cite{pollington_velani_2002}, and others, this result was generalized to badly
approximable points in $\RR^N$, $N>1$. In particular, the result
in~\cite{pollington_velani_2002} states that the set
$$
\bad(i_1,\ldots, i_N):=\left\{(x_1,\ldots,x_N)\in \RR^N:\;\exists \, c>0,
\max\{||qx_1||^{\frac{1}{i_1}},\ldots ||qx_N||^{\frac 1 {i_N}}\}\ge \frac c q \; \forall
q\in \NN\right\}
$$
has full Hausdorff dimension $N$, where $i_1,\ldots,i_N$ are any
strictly positive real numbers satisfying $i_1+\ldots+i_N=1$. Here,
$|| \, . \, ||$ denotes the distance to the nearest integer.
Finally, in~\cite{kristensen_thorn_velani_2006} a quite general
theory of badly approximable sets was constructed. It allows one to
establish full Hausdorff dimension results for a quite general class
of sets living in arbitrary compact metric spaces as long as certain
structural conditions are satisfied. The framework developed
encompasses  the results of both Jarn\'ik and Pollington \& Velani
as described above. Broadly speaking, the sets considered
in~\cite{kristensen_thorn_velani_2006} consist of points in a metric
space that avoid a given family of subsets of the metric space; that
is, points which cannot be easily approximated by this family of
subsets. Naturally, such sets were also referred to as badly
approximable. We give further details of this theory in
Section~\ref{sec_genbad}.


Recently, various sets were introduced within the theory of Diophantine approximation that on one hand could be naturally associated with the notion of badly approximable sets, but on the other hand do not seem to be covered
by the framework of~\cite{kristensen_thorn_velani_2006}. 
For example, in the landmark
paper~\cite{badziahin_pollington_velani_2011} the authors showed that
the set $\bad(i,j)$ intersected with any vertical line in $\RR^2$ is either empty
or has Hausdorff dimension equal to one. Whether this intersection is empty or has full dimension depends only upon a Diophantine property of the vertical line parameter. Later,
Beresnevich~\cite{beresnevich_2013} showed that for any
non-degenerate manifold $\MMM\subset \RR^N$,
$$
\dim(\bad(i_1,\ldots,i_N)\cap \MMM) = \dim \MMM,
$$
or in other words $\bad(i_1,\ldots,i_N)\cap \MMM$ has full Hausdorff
dimension.

One of the aims of this paper is to develop a framework in the
theory of badly approximable sets which will cover these new
results. In addition, we will show in detail (see Section
\ref{sec_app}) exactly how our theory can be used to attack various
problems in the field of Diophantine approximation, some of them old
and some of them previously open.


Our methodology appeals to the idea of \textit{generalised Cantor
sets in $\RR^N$}, which first appeared within the proofs
of~\cite{badziahin_pollington_velani_2011} and whose concept was
developed in the subsequent paper~\cite{badziahin_velani_2011}. The
construction of generalised Cantor sets has formed a basis for
establishing various difficult problems in the field of Diophantine
approximation. Many of these problems had proven resistant to
previous methods. For example, in \cite{badziahin_2013} generalised
Cantor sets were utilised to show that the set of points
$(x,y)\in\RR^2$ satisfying
\begin{equation}
\label{LC}
\liminf_{q\to \infty} q\cdot \log q\cdot \log\log q\cdot ||qx||\cdot
||qy||>0
\end{equation}
is of maximal Hausdorff dimension $2$ - a set not falling within the scope of \cite{kristensen_thorn_velani_2006}. This result represented significant progress in the investigation towards the famous Littlewood Conjecture, which states that the set of $(x,y)\in\RR^2$ satisfying \eqref{LC}, but with the `$\log q \cdot \log\log q$' term removed, is empty.

Whilst the Littlewood Conjecture is considered one of the most
profound an evasive problems in all of Diophantine approximation, in
recent years there has been much interest in a relatively new and
related problem. In 2004, de Mathan and Teuli\'{e} \cite{dMT}
proposed the following. Let $\mathcal{D}=(d_n)_{n\in \NN}$
be a sequence of positive integers greater or equal to 2 and let
$$
D_0:=1;\quad  D_n:=\prod_{k=1}^n d_k.
$$
Then, define the `pseudo-norm' function
$|\cdot|_\mathcal{D}:\NN\rightarrow \RR_{\ge 0}$ by
\[
|q|_\mathcal{D} = \min\{ D_n^{-1} : q\in  D_n\mathbb{Z} \}.
\]
If $\mathcal{D}=(p)_{n\in\NN}$ is a constant sequence for some
prime number $p$ then $|\cdot|_\mathcal{D}=|\cdot|_p$ is the usual
$p-$adic norm. The de Mathan-Teuli\'{e} Conjecture, often referred
to as the `Mixed' Littlewood Conjecture, is the assertion that for
any sequence $\mathcal{D}$ and for every $x\in\RR$ we have
\begin{equation}\label{MLC}
\liminf_{q\to\infty} q\cdot|q|_\mathcal{D}\cdot||q x||  \: = \: 0.
\end{equation}

In~\cite{badziahin_velani_2011}  generalised Cantor sets were utilised to show that the set
of real numbers $x\in\RR$ satisfying
\begin{equation}\label{MLClogs}
\liminf_{q\to \infty} q\cdot \log q\cdot \log\log q\cdot |q|_\mathcal{D}\cdot
||qx||>0
\end{equation}
is of maximal Hausdorff dimension and represents the state of the
art in results of this type (for arbitrary $\mathcal{D}$). An
application of the framework developed in our paper shows that for
sequences $\DDD$ growing sufficiently quickly
statement~\eqref{MLClogs} can be significantly improved. In
particular, we show that for any monotonic function $g\;:\NN\;\to
\RR_{\ge 0}$ tending to infinity and every sequence $\DDD =
(d_i)_{i\in\NN}$ such that
$$
\lim_{i\to\infty} \frac{g(d_{i+1})}{g(d_i)}= \infty,
$$
the set of real numbers $x\in\RR$ satisfying
\begin{equation}\label{MLCg}
\liminf_{q\to \infty} q\cdot g(q)\cdot |q|_\mathcal{D}\cdot
||qx||>0
\end{equation}
is of maximal Hausdorff dimension. Note that we may choose $g$
tending to infinity as slowly as we wish. The only previously known
result of this type `beating' the rate of approximation in
\eqref{MLClogs} was proved in \cite{badziahin_velani_2011}, where it
was shown that the set of $x\in\RR$ satisfying \eqref{MLCg} with
$g(q)=\log\log q \cdot \log\log\log q$ has maximal Hausdorff
dimension in the specific case
$\mathcal{D}=\{2^{2^n}\}_{n=0}^\infty$.

Our paper also extends the concept of generalised Cantor sets to the
setting of general metric spaces. This allows us to utilise modern
techniques in setups to which they did not previously extend. As
one of the applications considered in this paper we consider the
space~$\ZZ_p$ of $p$-adic integers. For $N \in \NN$, the set $\bad_p(N)$ of so-called
\textit{badly approximable $p$-adic vectors} is defined as the collection of points $(x_1,\ldots,x_N)\in \ZZ_p^N$ for which there
exists a constant $c>0$ satisfying
\begin{equation}\label{def_badp}
\max\{|qx_1-r_1|_p,\ldots |qx_N-r_N|_p\}\ge c \cdot \max\{ |r_1|,
\ldots, |r_N|, |q| \}^{-\frac{N+1}{N}}
\end{equation}
for every $(r_1, \cdots, r_N, q)\in \ZZ^N \times \ZZ\setminus \{0\}$.
The set $\bad_p(1)$ was shown to have maximal Hausdorff dimension by Abercrombie~\cite{Aber} in 1995. In 2006, this result was extended using the broad framework of~\cite{kristensen_thorn_velani_2006}, where it was shown that the set $\bad_p(N)$ has maximal Hausdorff dimension~$N$. However, nothing is known about any
stronger properties of $\bad_p(N)$, such as whether it is winning with respect to Schmidt's game, which we now introduce. We show that at
the very least $\bad_p(N)$ satisfies the amenable properties enjoyed by the
`winning sets' occurring in Schmidt's game. Establishing these
properties had previously appeared out of reach.

%

We also find answers to some other new problems from the field of
Diophantine approximation, such as questions relating to the $\times 2,
\times 3$ problem, and questions relating to the behaviour of the Lagrange
constant of multiples of a given real number as posed in \cite{bugeaud_2015}.

\subsection{Winning sets and countable intersections}\label{sec:winning}
Given a ball $B$ we write $\rad(B)$ and $\diam(B)$ for the radius
and the diameter of $B$ respectively. By $\cent(B)$ we denote the
center of $B$.

Another remarkable property of the set $\bad$ was discovered by
Schmidt in a series of works finalised in~\cite{schmidt_1966}. It
can be described in terms of Schmidt's so called  \textit{$(\alpha,\beta)$-game}.
Suppose two players Alice (A) and Bob (B) play the following game with two fixed real parameters $0<\alpha,\beta<1$.
Bob starts by choosing an arbitrary closed ball $B_1\subset\RR^N$. Then
Alice and Bob take turns in choosing closed
balls in a nested sequence (Alice chooses balls $A_i$ and Bob chooses balls $B_i$),
$$
B_1\supset A_1\supset B_2 \supset A_2\supset\cdots,
$$
whose radii satisfy
$$
\rad(A_i) = \alpha \cdot \rad(B_i),\quad \rad(B_{i+1}) = \beta \cdot \rad(A_i),\quad \forall
i\in \NN.
$$
A set $E\subset \RR^N$ is called \textit{$(\alpha,\beta)$-winning} if Alice
has a strategy that ensures that
$$
\bigcap_{i=1}^\infty B_i = \vp  \, \in E.
$$
Finally, we say that $E\subset\RR^N$ is\textit{ $\alpha$-winning} if it is
$(\alpha,\beta)$-winning for all $0<\beta<1$ and \textit{winning } if it is $\alpha$-winning for some $\alpha \in (0,1)$.

Surprisingly, Schmidt was able to show that the set $\bad$ is winning as a subset of $\RR$, and further, that  all winning sets in Euclidean space satisfy some remarkable properties:
\begin{itemize}
    \item[(W1)] Any winning set is dense and has full Hausdorff dimension.
    \item[(W2)] Any countable intersection of $\alpha$-winning sets is $\alpha$-winning.
    \item[(W3)] The image of any winning set under a bi-Lipschitz map is again
    winning.
\end{itemize}

Many other sets including $\bad(1/N,1/N,\ldots,1/N)$ (which for brevity we
refer as $\bad_N$) have been proven to be winning \cite{schmidt_1980}. Most recently,  in an exceptional paper by An \cite{an_2013} it was shown that the set $\bad(i,j)$ is winning. This provided a second proof of the Badziahin-Pollington-Velani Theorem, the main result of~\cite{badziahin_pollington_velani_2011}, which in turn established a long standing conjecture of Schmidt. It appears that many sets that according
to~\cite{kristensen_thorn_velani_2006} fall into the category of
badly approximable sets are indeed winning.

Several variations of Schmidt's $(\alpha,\beta)$-game have been suggested who's analogous winning sets still satisfy the properties (W1)~--~(W3) of
classical winning sets. On the further development of the subject we
refer the reader to~\cite{kleinbock_weiss_2010, mcmullen_2010, schmidt_1980} and the references therein, and to Section \ref{sec:relationship} of this paper for a partial overview. One disadvantage of working with topological games of this type is that it is often quite difficult to prove a set is winning. A major aim of this paper is to develop a variation of the category of winning sets such that properties (W1)~--~(W3) are still satisfied for sets in this new category and the conditions for inclusion this category are rather easier to check. In particular, we define so called \textit{Cantor-winning} sets in Section~\ref{sec_winning}; a category of sets each of whom contain a generous supply of generalised Cantor sets.

As examples, the set consisting of the points in $\bad(i,j)$ lying
on certain vertical lines which appeared
in~\cite{badziahin_pollington_velani_2011} turns out to be
Cantor-winning, and so does the set of points in
$\bad(i_1,\ldots,i_N)$ lying on non-degenerate curves as described
in~\cite{beresnevich_2013}. Therefore, these sets each satisfy
properties (W1)~--~(W3). We discuss these results in more details in
Section~\ref{sec_app}, along with some other far reaching
applications.

\subsection{The idea of generalised Cantor sets}

The basic premise for the construction of generalised Cantor sets in an arbitrary metric space is the standard middle-third Cantor set
construction. We now describe this process and then discuss what requirements should be satisfied in order to generalize the construction to an arbitrary metric space. The classical Cantor set is realised as follows. We start
with the unit interval $I_0 = [0,1]$. The first step of the process is to split the interval $I_0$ into three
intervals of equal length and remove the open middle interval. This leaves a union $I_1 = [0,1/3]\cup [1/3,1]$ of two disjoint closed intervals which survive the first step. We
recursively repeat this procedure for each of the remaining
intervals, each time removing the open middle third from every interval in the union, to produce a sequence $(I_2,I_3,\ldots)$ of sets. Each $I_n$ will consist precisely of the disjoint union of the $2^i$ closed intervals that survive the $i$-th step of the removal procedure. The
classical middle-third Cantor set $\KKK$ is then defined as
$$
\KKK:= \bigcap_{n=0}^\infty I_n.
$$
The set  $\KKK$ is well known to be uncountable and have Hausdorff dimension $\log2/\log3$. Surely for $I_0$ we can take any interval instead of
$[0,1]$ and the final set $\KKK$ will still satisfy the same properties.

In an arbitrary metric space $X$ the (metric) balls will play the role of
intervals in $\RR$. One needs to define the rules describing how each
surviving ball should be split into smaller pieces in the next step of the construction. For example, when $X=\RR$
we may generalise the set $\KKK$ by splitting intervals into $R$ closed pieces of equal length at each step for some $R \geq 3$, or even varying the number of intervals created during each step of the procedure. In $\RR^N$ we can
take square boxes (that is, balls in the sup-norm metrics) and split
them into $R^N$ smaller boxes. For a general metric
space~$X$ we will need to describe how every ball is split into smaller
balls. In order for such a process to be meaningful (or even possible) we must enforce some kind of structure on $X$ which allows for such a splitting procedure to take place. In Section~\ref{sec1} we define an extremely general class of metric spaces possessing such a \textit{splitting structure}.

Returning to $\KKK$ for a moment, we may also generalise its
construction by varying the number of intervals removed at each step. However, this should be done with
care. For example, if in a classical Cantor set construction instead
of one interval we remove two of them on each step (let's say a
middle and left one) then every step will leave just one
interval: $I_1 = [2/3,1]$ , $I_2 = [8/9,1]$ and so on. In this case
$\bigcap_{n=0}^\infty I_n$ is a single point, which is probably not as
interesting as $\KKK$. For this reason we need to control the number of
 intervals produced in the splitting process together with the number of removed intervals in each step in
order to get a non-trivial generalised Cantor set at the end. We
provide reasonable restrictions on these numbers in
Section~\ref{sec_gencant}, although they are essentially the same as
in~\cite{badziahin_velani_2011}. The key point of the described
procedure is that we shall allow the number of intervals removed at
each step of the removal process to depend on the entirety of the
construction thus far, not just upon the specific step as in the
classical Cantor set construction.

\section{Splitting structure and metric spaces}\label{sec1}

We now describe sufficient conditions on a metric space $X$ for a generalised Cantor construction to be possible. Denote by  $\BBB(X)$ the
set of all closed (metric) balls in $X$.

We define a \textit{splitting structure} on a metric space $X$ (with metric $\vd$) as a quadruple~$(X,\SSS, U,f)$, where
\begin{itemize}
\item $U\subset \NN$ is an infinite multiplicatively closed set such that if
$u,v\in U$ and $u\mid v$ then $v/u\in U$;
\item $f\;:\; U\to \NN$ is an absolutely multiplicative arithmetic function;
\item $\SSS\;:\BBB(X)\times U \to\BBB(X)$ is a map defined in such a way that for every ball $B\in\BBB(X)$ and
$u\in U$, the set $\SSS(B,u)$ consists solely of balls $b_i\subset B$ of
radius $\rad(B)/u$.
\end{itemize}
Additionally, we require all these objects to be linked by the
following properties
\begin{enumerate}
\item[(S1)] $\#\SSS(B,u) = f(u)$;
\item[(S2)] If $b_1,b_2\in \SSS(B,u)$ and $b_1\neq b_2$ then $b_1$ and $b_2$ may only intersect on the
boundary; i.e.
$$
\vd(\cent(b_1),\cent(b_2)) \ge \frac{2\cdot \rad(B)}{u};
$$
\item[(S3)] For all $u,v\in U$,
$$
\SSS(B,uv) = \bigcup_{b \,\in \,\SSS(B,u)} \SSS(b,v).
$$
\end{enumerate}

\noindent{\bf Remark.} Not all metric spaces possess a splitting
structure. For example, it is easy to check that if $f\not\equiv 1$
then $X$ must be infinite. On the other hand in the case $f \equiv
1$ we always have that $\SSS(B, u)$ consists of one ball. This case
is not very interesting and we call such a splitting structure {\it
trivial}. Furthermore, given a metric space  $X$, there usually
exist some restrictions on the growth of $f$ for the splitting
structure $(X,\SSS,U,f)$ to exist. For example, when $X=\RR^N$
properties (S1) and (S2) imply that we must have $f(u)\le u^N$.

Note also that $\SSS(B,u)$ does not necessarily form a cover of $B$. However in the
cases we are mostly interested in this property will be satisfied.

\subsection{Some examples}
\begin{itemize}
\item[(a)] Let $X=\RR^N$ with $\vd(\vx,\vy):=|\vx-\vy|_\infty$, $U =
\NN$, $f(u) = u^N$ and $\SSS(B,u)$ be defined as follows: $B$ is cut
into $u^N$ equal square boxes which edges have length $u$ times less
than the edges length of $B$. One can easily check that
$(\RR^N,\SSS, \NN, f)$ satisfies properties (S1)~--~(S3). We call this the
\textit{canonical }splitting structure for $\RR^N$.

\item[(b)] Let $X = \QQ_p^N$ with $\vd(\vx,\vy):=
\max_{1\le i\le N}\{|x_i-y_i|_p\}$, $U = \{p^k\;:\; k\in\ZZ_{\ge
0}\}$, $f(p^k) = p^{Nk}$ and $\SSS(B,p^k)$ be defined as the set of
all disjoint balls in $B$ of radius $\rad(B)/p^k$. Again properties
(S1) -- (S3) are easily verified, so $(\QQ_p^N, \SSS, U, f)$ is a
splitting structure. We call this splitting structure
\textit{canonical} for $\QQ_p^N$.

\item[(c)] We give one more exotic example. Let $X = \RR$, $U = \{3^k\;:\;
k\in\ZZ_{\ge 0}\}$, $f(3^k) = 2^k$. Define $\SSS(B,3)$ as follows: we divide the interval $B$ into 3 pieces of equal length
and remove the open third in the middle. $\SSS(B,3^k)$ for $k>1$ is constructed
inductively with help of property~(S3). It is easily verified that $(\RR, \SSS, U, f)$ indeed forms a splitting structure.
\end{itemize}

We will refer to these examples throughout the paper.

\subsection{The set $A_\infty(B)$}

A splitting structure on a metric space naturally exhibits a Cantor-like structure. For $u\in U$ and $B\in\BBB(X)$ define the set
$$
A_u(B):= \bigcup_{b\, \in\, \SSS(B,u)} b.
$$
By property (S2), if $u,v\in U$ and $u\mid v$ then $A_v(B) \subset
A_u(B)$. Moreover, if $X$ is complete then for every
sequence $(u_i)_{i\in\NN}$ with $u_i \in U$ such that $u_i\mid u_{i+1}$ the set
$\bigcap_{i=1}^\infty A_{u_i}(B)$ is non-empty. Moreover, it has
non-empty intersection with each ball from $\SSS(B,u_i)$ and the following property also holds.

\begin{theorem}\label{th1}
Let  $(X, \SSS, U, f)$ be a splitting structure. Then for every two sequences
$(u_i)_{i\in\NN}$ and $(v_i)_{i\in\NN}$ with $u_i, v_i\in U$ such that
$u_i\mid u_{i+1}, v_i\mid v_{i+1}$ and
$u_i,v_i\stackrel{i\to\infty}{\longrightarrow} \infty$ one has
$$
\bigcap_{i=1}^\infty A_{u_i}(B)  = \bigcap_{i=1}^\infty A_{v_i}(B).
$$
\end{theorem}
\proof Suppose that the contrary is true; that is, suppose there exists $x\in X$ such that
$$
x\in \bigcap_{i=1}^\infty A_{u_i}(B),\quad\mbox{but}\quad x\not\in
\bigcap_{i=1}^\infty A_{v_i}(B).
$$
Since $\bigcap_{i=1}^\infty A_{v_i}(B)$ is closed, there exists a real number $\delta>0$ and a ball
$B(x,\delta) \in \BBB(X)$ around $x$ such that
$$
B(x,\delta)\cap \bigcap_{i=1}^\infty A_{v_i}(B) = \emptyset.
$$

By the construction of the sets $A_{u_i}(B)$ there exists $u_i$ and a ball
$b\in \SSS(B, u_i)$ of radius less than $\delta/2$ such that $x\in
b$. Therefore, $b\subset B(x,\delta)$. Finally, by property~(S2), $A_{u_i\cdot v_j}(B)$ has non-empty intersection with $b$ and is a subset of $A_{v_j}(B)$ for every
$j\in\NN$. Whence, by taking intersection we get
$$
\emptyset \neq b\cap \bigcap_{j=1}^\infty A_{u_i\cdot v_j}(B)\subset
B(x,\delta)\cap \bigcap_{j=1}^\infty A_{v_j}(B)
$$
and reach a contradiction.\endproof

The crux of Theorem \ref{th1} is that an infinite intersection
$\bigcap_{i=1}^\infty A_{u_i}(B)$ depends only on the ball $B$ and the
splitting structure on $X$, but not on the particular sequence
$(u_i)$. We denote this intersection by $A_\infty(B)$. In further
discussion it will always be assumed that $X$ is a complete metric
space and so the notion $A_\infty(B)$ will be always correctly
defined.

It is readily observed that for any trivial splitting structure
$A_\infty (B)$ consists of just a single point. Also, one can easily
check that for the canonical splitting structure on $\RR^N$ and for
the canonical splitting structure on $\QQ_p^N$ (examples (a) and
(b)), we have that $A_\infty(B) = B$. In the case of more exotic
splitting structure from example (c) one can check that
$A_\infty(B)$ is a standard middle-third Cantor set $\KKK(B)$ whose
construction starts with interval $I_0=B$. As previously discussed,
the set $\KKK(B)$  is compact. Indeed, the set  $A_\infty(B)$ is
compact for each of the examples (a)~-~(c). We now demonstrate that
this property is actually ubiquitous.

\begin{theorem}\label{th2}
Let   $(X, \SSS, U, f)$ be a splitting structure. Then, for any ball
$B \in \BBB(X) $ the set~$A_\infty(B)$ is compact.
\end{theorem}

\proof For a trivial splitting structure the result is obvious.
Therefore, assume that the splitting structure is non-trivial. Fix
parameter $v\in U$ with $v>1$. Consider a cover $\bigcup_\alpha O_\alpha$
of $A_\infty (B)$ by open sets $O_\alpha$. To each $O_\alpha$ we may
associate a subset of balls from $\bigcup_{i=1}^\infty\SSS(B,v^i)$
such that every ball $b$ from this subset lies entirely inside
$O_\alpha$. In particular, let
$$
\DDD(O_\alpha):= \left\{b\in \bigcup_{i \in \NN} \SSS(B,v^i)\;:\;
b\subset O_\alpha\right\}.
$$
Obviously, if $b\in\DDD(O_\alpha)$ then every ball
$b'\in\SSS(b,v^j)$ for $j \in \NN$ is also in $\DDD(O_\alpha)$.

If for some $i\in\NN$ every ball from $\SSS(B,v^i)$ is in one of the
sets $\DDD(O_\alpha)$ (for some $\alpha$) then there is a finite
subcover of $A_\infty(B)$. Indeed, for every ball $b$ from the
finite set $\SSS(B,v^i)$ we associate one element $O_\alpha$ from
the cover such that $b\in\DDD(O_\alpha)$. Assume now that this is
not the case. Then there exists a sequence  $(b_i)_{i\in
\NN}$ of balls such that $b_i\in\SSS(B,v^i)$, $b_i\supset b_{i+1}$ such that
none of these balls are in $\DDD(O_\alpha)$ for any $\alpha$. Since
$X$ is complete we have $\bigcap_{i=1}^\infty b_i=x$ is a single
point. It must be covered by one of the open sets $O_\alpha$ and so
there must exist $\epsilon>0$ such that $B(x,\epsilon)\subset
O_\alpha$. Moreover, for $i$ large enough we must have $b_i\subset
B(x,\epsilon)\subset O_\alpha$, a contradiction.
\endproof

As in the case of the middle-third Cantor set it is desirable to determine the Hausdorff dimension of the set
$A_\infty(B)$. In order to compute this in general we require the metric space $X$ to satisfy one further condition.

\begin{itemize}
\item[(S4)] There exists an absolute constant $C(X)$ such that any
ball $B \in \BBB(X)$ cannot intersect more than $C(X)$ disjoint open balls of
the same radius as $B$.
\end{itemize}

It is easy to check that the metric spaces $\RR^N$ and $\QQ_p^N$ satisfy Condition
(S4), and so all of the examples (a)~-~(c) satisfy it. This condition is sufficient to precisely calculate the Hausdorff dimension of $A_\infty(B)$.

\begin{theorem}\label{th3}
    Let  $(X, \SSS, U, f)$ be a splitting structure. Then for any $B \in \BBB(X)$ we have
$$
\dim A_\infty (B) \le \liminf_{u\to\infty; u\in U} \frac{\log
f(u)}{\log u}.
$$
Moreover, if $X$ satisfies condition (S4) then $\log f(u) / \log u$
must be a constant and
$$
\dim A_\infty (B) = \frac{\log f(u)}{\log u}.
$$
\end{theorem}

\proof Determining the upper bound for the Hausdorff dimension is relatively easy. We may simply consider the trivial cover $\SSS(B,u)$ of $A_\infty (B)$. We have
$$
\sum_{b\in\SSS(B,u)} (\rad(b))^d \asymp f(u)\cdot u^{-d},
$$
and so for $d>\liminf \log f(u)/ \log u$ one can find a sequence of integers $u_i\in
U$ and $\epsilon>0$ such that
$$
\frac{\log f(u_i)}{\log u_i}<d-\epsilon.
$$
Therefore,
$$
f(u_i)\cdot u_i^{-d} < f(u_i)\cdot u_i^{-\frac{\log f(u_i)}{\log
u_i}-\epsilon} = u_i^{-\epsilon}\; \stackrel{i\to\infty}
\longrightarrow\; 0
$$
and this gives us a required upper bound on $\dim A_\infty (B)$.

The inverse inequality requires a bit more effort. Consider $u\in U$ and
let $d = \log f(u) / \log u$. If we prove that
\begin{equation}\label{eq1}
\inf \left\{\sum_i (\rad(B_i))^d \;:\; \bigcup_i B_i\mbox{ is a
cover of }B\right\}>0
\end{equation}
then we will have that $\dim A_\infty(B)\ge d$ as required.

Let $\bigcup_\alpha B_\alpha$ be an arbitrary cover of $A_\infty(B)$ by
open balls $B_\alpha$. Then, since $A_\infty(B)$ is compact one can choose a
finite subcover $\bigcup_{i=1}^n B_i$. This procedure only
decreases the value on the left hand side of~\eqref{eq1} and so if we can show that
\begin{equation}\label{eq1finite}
\inf \left\{\sum_i (\rad(B_i))^d \;:\; \bigcup_i B_i\mbox{ is a
    finite cover of }B\right\}>0
\end{equation}
then we are done. Given a finite subcover $\bigcup_{i=1}^n B_i$ one may without loss of generality
assume that $\rad(B_i)\le \rad(B)$ for each $1\le i\le n$, for otherwise~\eqref{eq1finite} has an obvious positive infimum equal to
$(\rad(B))^d$.

Consider an individual element $B_i$ of the subcover. Let $m_i\in\ZZ_{\ge 0}$ take a value such that
$$
\frac{\rad(B)}{u^{m_i+1}}< \rad(B_i)\le \frac{\rad(B)}{u^{m_i}}.
$$
Then, take all balls $B_{i,1}, B_{i,2},\ldots, B_{i,s_i}$ from the collection $ \SSS(B,
u^{m_i})$ which have non-empty intersection with $B_i$. By
condition~(S4) we must have $s_i\le C(X)$ and so
$$
(\rad(B_i))^d\ge \frac{1}{C(X)\cdot u^d}\sum_{s=1}^{s_i}
(\rad(B_{i,s_i}))^d.
$$
Replacing each $B_i$ by $B_{i,1},\ldots,B_{i,s_i}$ one can easily
check that
\begin{equation}\label{eq2}
\bigcup_{i=1}^n \bigcup_{s=1}^{s_i} B_{i,s}
\end{equation}
is still a cover of $A_\infty(B)$ and it solely consists of balls
from $\bigcup_{i=1}^\infty \SSS(B,u^i)$. The value in~\eqref{eq1finite} does not decrease more than $C(X)\cdot u^d$ times compared to the
initial cover $\bigcup_{i=1}^n B_i$.

Now notice that by the definition of $d$, for every ball $B'$ one has
$$
(\rad(B'))^d = \sum_{b\in\SSS(B',u)}(\rad(b))^d.
$$

In other words, if in a cover one replaces one ball $B'$ by all balls
in $\SSS(B',u)$, the value~\eqref{eq1} does not change. We use this
observation and replace if necessary every ball $B_{i,s}$
from~\eqref{eq2} by balls from $\SSS(B_{s,i}, u^{k_{s,i}})$ for some
$k_{s,i}$ to guarantee that all balls in the resulting cover are
from $\SSS(B, u^k)$ for some fixed $k\in\NN$. Since it is still a
cover of $B$, all these balls together comprise the whole set
$\SSS(B,u^k)$. Finally,
$$
\sum_{i=1}^n (\rad(B_i))^d \ge \frac{1}{C(X)\cdot u^d}\cdot
\sum_{b\in \SSS(B,u^k)}(\rad(b))^d = \frac{(\rad(B))^d}{C(X)\cdot
u^d}> 0.
$$
The claim is achieved, therefore we get
$$
\dim A_\infty(B) \le \frac{\log f(u)}{\log u}.
$$
To finish the proof we take an arbitrary $u\in U$ and combine the
last statement with the lower bound for $\dim A_\infty(B)$ we got
before. \endproof

\begin{corol}\label{corol1}
If $X$ satisfies condition~(S4) and $U$ contains at least two
multiplicatively independent numbers then $f(u)$ must be of the
form: $f(u) = u^d$ where $d\in \QQ$, $d>0$.
\end{corol}

\proof By Theorem~\ref{th3},
$$
d = \frac{\log f(u)}{\log u}
$$
is a constant. Therefore $f(u) = u^d$. If $u_1,u_2$ are two
multiplicatively independent elements of $U$, then $u_1^d$ and
$u_2^d$ can both be integer only if $d\in\QQ$. Finally $d>0$ since
otherwise $u^d\not\in\ZZ$. \endproof

Notice that if $U$ is generated by one positive integer number
$u_0$ then for every $u=u_0^n$ we have $f(u) = f(u_0)^n = u^d$ where
$d = \frac{\log f(u_0)}{\log u_0}$. Again $f$ is of the form $f(u) =
u^d$, however in this case we do not necessarily have that
$d\in\QQ$.

\section{Generalised Cantor sets}\label{sec_gencant}

Let $(X,\SSS,U,f)$ be a splitting structure on $X$. We now introduce the precise definition of generalised Cantor sets
 in the context of  this splitting structure, for which we appeal heavily to the ideas presented in \cite{badziahin_velani_2011}.

Fix some closed ball $B \in \BBB(X)$, let
$$
\RRR: = (R_n)_{n\in\ZZ_{\ge 0}},\; R_n\in U
$$
be a sequence of natural numbers and let
$$
\vr:= (r_{m,n}),\; m,n\in\ZZ_{\ge 0},\; m\le n
$$
be a two parameter sequence of real numbers.

\noindent{\bf Construction.} We start by considering the set
$\SSS(B,R_0)$. The first step in the construction of a generalised
Cantor set involves the removal of at most $r_{0,0}$ balls $b$ from
$\SSS(B,R_0)$. We call the resulting set $\BBB_1$. Balls in $\BBB_1$
will be referred as (level one) survivors. Note that we do not
specify the removed balls, just give an upper bound for their
number. For consistency we also define $\BBB_0:=\{B\}$.

In general, for $n\ge 0$, given a collection $\BBB_n$ we construct a
nested collection $\BBB_{n+1}$ using the following two operations:
\begin{itemize}
\item Splitting procedure: Compute the collection of candidate balls
$$
\III_{n+1}:= \bigcup_{B_n\in\BBB_n}\SSS(B_n, R_n).
$$
\item Removing procedure: For each ball $B_n\in \BBB_n$ we remove at
most $r_{n,n}$ balls $B_{n+1}\in\SSS(B_n,R_n)$ from $\III_{n+1}$. Let
$\III_{n+1}^n\subseteq \III_{n+1}$ be the collection of balls that
remain. Next, for each ball $B_{n-1}\in \BBB_{n-1}$ we remove at
most $r_{n-1,n}$ balls $B_{n+1}\in \SSS(B_{n-1},R_nR_{n-1})\cap
\III_{n+1}^n$. Let $\III_{n+1}^{n-1}$ be the collection of balls
that remain. In general for each $B_{n-k}\in\BBB_{n-k}$ ($1\le k \le
n$) we remove at most $r_{n-k,n}$ balls $B_{n+1}\in
\SSS(B_{n-k},\prod_{i=0}^k R_{n-i})\cap \III_{n+1}^{n-k+1}$ and
define $\III_{n+1}^{n-k}\subseteq \III_{n+1}^{n-k+1}$ to be the
collection of balls that remain. Finally, $\BBB_{n+1}:=
\III_{n+1}^0$ then becomes the desired collection of (level $n+1$) survivors.
\end{itemize}

The two operations above allow us to construct a nested sequence of  collections $\BBB_n$ of closed
balls. Consider the limit set
$$
\KKK(B,\RRR, \vr):= \bigcap_{i=1}^\infty \bigcup_{b\in\BBB_n} b.
$$
The set $\KKK(B,\RRR,\vr)$ will be referred to as a \textit{generalised $(B,\RRR,\vr)$-Cantor set} on $X$.

Note that the triple $(B,\RRR,\vr)$ does not uniquely determine
$\KKK(B,\RRR,\vr)$. There is a large degree of freedom in the choice of balls $B_{n+1}$ removed in the construction procedure. Consequently, one can look at the property of being a generalised $(B,\RRR,\vr)$-Cantor set as a property of the set $\KKK\subset X$, rather than as a self contained definition: we say a set $\KKK$ is
\textit{a generalised Cantor set} if it can be constructed by the procedure described above for some triple $(B,\RRR,\vr)$. In this case, we may refer to $\KKK$ as being $(B,\RRR,\vr)-Cantor$ if we wish to make such a triple explicit and write $\KKK=\KKK(B,\RRR,\vr)$.

\subsection{Properties of $\KKK(B,\RRR,\vr)$}

Generalized $(B,\RRR,\vr)$-Cantor sets in any complete metric space $X$ satisfy the same
desirable properties as proved in~\cite{badziahin_velani_2011}. Furthermore, many of the proofs translate from the Euclidean setting to the case of arbitrary metric spaces with only slight modification. We now exhibit these properties,  but will provide the proof only if it significantly differs from the analogous methods outlined in~\cite{badziahin_velani_2011}.

\begin{theorem}[See Theorem 3 in~\cite{badziahin_velani_2011}]\label{th_cantor1}
Given a generalised Cantor set $\KKK (B,\RRR,\vr) $ in a complete metric space $X$, let
\begin{equation}\label{def_t0}
t_0:=f(R_0)-r_{0,0}
\end{equation}
and for $n \ge 1 $ let
\begin{equation}\label{def_tn}
t_n:=f(R_n)-r_{n,n}-\sum_{k=1}^n \frac{r_{n-k,n}}{\prod_{i=1}^k
t_{n-i}}    \,    .
\end{equation}
Suppose that  $t_n>0$ for all $n\in\ZZ_{\ge 0}$. Then,
$$
\KKK ({\rm B},\RRR,\vr)  \neq \emptyset \ .
$$
\end{theorem}

\begin{theorem}[See Theorem 4 in~\cite{badziahin_velani_2011}]\label{th_cantor2}
Let a complete metric space $X$ satisfy condition~(S4). Given a generalised Cantor set $\KKK ({\rm
B},\RRR,\vr) \subset X $, suppose that the parameters $\RRR$ and $\vr$ satisfy
the following conditions:
\begin{itemize}
\item $f(R_n)\ge 4$ for all $n\in\ZZ_{\ge 0}$;
\item for every $\delta>0$ there exists $n(\delta)$ such that for
every $n>n(\delta)$,
\begin{equation}\label{ineq_lem2}
\prod_{i=0}^n R_i^\delta >R_n^s,
\end{equation}
where $s = \liminf\limits_{n\to\infty} (\dim A_\infty(B) -
\log_{R_n}2)$;
\item For every $n\in\ZZ_{\ge 0}$,
\begin{equation}\label{cond_th2}
\sum_{k=0}^n \left(r_{n-k,n}\prod_{i=1}^k
\left(\frac{4}{f(R_{n-i})}\right)\right)\le \frac{f(R_n)}{4}.
\end{equation}
\end{itemize}
Then
$$
\dim \KKK ({\rm B},\RRR,\vr)  \ge s.
$$
\end{theorem}

\noindent{\bf Remark.} It is unclear to the authors as to whether
condition~\eqref{ineq_lem2} is absolutely necessary. For example, in the
corresponding Theorem~4 from~\cite{badziahin_velani_2011} this
condition is not needed. On the other hand it may be the property of
the canonical splitting structure of $\RR$ that
makes~\eqref{ineq_lem2} superfluous. Whilst the proof of  Theorem \ref{th_cantor2} is very similar to that in~\cite{badziahin_velani_2011}, we consider it to be quite
important, especially with regard to the need for condition~\eqref{ineq_lem2}.  For this reason, and for the sake of completeness, we briefly outline its
proof here.

Prior the proof we give a definition of \textit{local Cantor sets}, which provide the means by which one can prove most of the results in this section. A generalised Cantor
set $\KKK({\rm B},\RRR,\vr)$ is said to be \textit{local} if $r_{m,n}=0$
whenever $m\neq n$. Furthermore, we write $\LKK({\rm B},\RRR,\vs)$
for $\KKK({\rm B},\RRR,\vr)$ where
$$
\vs:=(s_n)_{n\in \ZZ_{\ge 0}}  \quad   {\rm and } \quad
s_n:=r_{n,n}.
$$

We will also need the following version of the mass distribution principle for general metric spaces $X$, a powerful tool for calculating lower bounds for Hausdorff dimension.

\begin{MDP}
Let $ \mu $ be a probability measure supported on a subset $E$ of a
metric space $X$. Suppose there are  positive constants $a, s $ and
$l_0$ such that
\begin{equation}\label{mdp_eq1}
\mu  ( B ) \le \, a \;   \diam(B)^s \; ,
\end{equation}
for any closed set $B$ with $\diam(B)\le l_0$. Then,  $\dim E \ge
s$.
\end{MDP}

One can check that it is sufficient to verify
property~\eqref{mdp_eq1} for all balls $B\in\BBB(X)$. Indeed, assume
that it is satisfied for balls. Consider an arbitrary set $S\subset X$
of diameter $\diam(S)\le l_0/2$. It is covered by a ball $B$ with
$\rad(B)\le \diam(S)$, so $\diam(B)\le l_0$. Then we have
$$
\mu(S)\le \mu(B)\le a\cdot \diam(B)^s \le a\cdot 2^s\cdot
\diam(S)^s.
$$
Therefore Property~\eqref{mdp_eq1} is satisfied then for an
arbitrary set $S$ with parameters $a':=a\cdot 2^s,$ $s':=s$ and $l_0':= l_0/2$. It follows that $\dim E\ge s'=s$.

The final prerequisites for the proof of Theorem~\ref{th_cantor2}
are a lower bound for the Hausdorff dimension of local Cantor sets
and a proof that certain generalised Cantor sets contain
sufficiently permeating local Cantor sets.

\begin{lemma}\label{lem_local}
Given   $\LKK({\rm B},\RRR,\vs)$, suppose that
$$
t_n:=f(R_n)-s_n   > 0  \quad \forall \  n\in\ZZ_{\ge 0}  \, .
$$
Furthermore, suppose the values $s_n$ and $R_n$ satisfy the
following conditions: there are positive constants $s$ and $n_0 $
such that for all $n>n_0$
\begin{equation}\label{ineq_lemloc}
R_n^s\le t_n   \,
\end{equation}
and for every $\delta>0$ there exists $n(\delta)>0$ such that
inequality~\eqref{ineq_lem2} is satisfied. Then
$$
\dim \LKK(B,\RRR,\vs) \ge s.
$$
\end{lemma}

\proof We construct  a  probability measure $\mu$ supported
on $\LKK(B,\RRR,\vs)$ in the standard manner. For any $B_n \in
\BBB_n$, we attach a weight $\mu(B_n)$ defined recursively as
follows.

\noindent For  $n=0$ let $$ \mu(B_0)\ := \frac{1}{\#\BBB_0}=1,  \ $$ and
for $n\ge 1$ define
\begin{equation}\label{beq4}
\mu(B_n)  \, :=  \, \frac{\mu(B_{n-1})}{\# \{B\in \BBB_n\;:\;
B\subset B_{n-1}\}}, \
\end{equation}
where $B_{n-1} \in \BBB_{n-1}$   is the unique ball such that
$B_n\subset B_{n-1}$. This procedure  inductively defines a mass on
any interval appearing in the construction of $\LKK(B,\RRR,\vs) $.
In fact, it can be easily demonstrated via induction that for every
$B_n\in \BBB_n$ we have
\begin{equation}\label{eq3}
\mu(B_n) \le \prod_{i=0}^{n-1} t_i^{-1}.
\end{equation}
This measure can be further extended to all Borel subsets of
$X$. We will call such a measure {\it a canonical measure} on
$\LKK(B,\RRR,\vs)$. It remains to show that $\mu$
satisfies~\eqref{mdp_eq1}. Consider an arbitrary ball $E$ of radius
not bigger than $\rad(B)$. Then there exists a positive integer
parameter $m$ such that
\begin{equation}\label{eq4}
\frac{\rad(B)}{\prod_{i=0}^m R_i} < \rad(E)\le
\frac{\rad(B)}{\prod_{i=0}^{m-1} R_i}.
\end{equation}

Now we estimate $\mu(E)$. First, notice that we have
$$
\mu(E)\le \sum_{b\in\BBB_m: \,\; b \, \cap \, E \, \neq \,  \emptyset} \mu(b).
$$
By Property~(S4) there are at most $C(X)$ balls $b\in\BBB_m$ such
that $b\cap E\neq \emptyset$. This, together with~\eqref{eq3}, gives us the upper bound
$$
\mu(E)\le C(X)\cdot \prod_{i=0}^{m-1} t_i^{-1}\le C(X)\cdot
\prod_{i=0}^{m-1} \frac{R_i^s}{t_i}\big/ \prod_{i=0}^{m-1}R_i^s
 \stackrel{\eqref{ineq_lem2}}{\le} C(X)\cdot
C_1(\delta)\frac{\rad(B)^{s-\delta}}{\prod_{i=0}^mR_i^{s-\delta}}\cdot
\prod_{i=0}^{m-1} \frac{R_i^s}{t_i},
$$
where
$$
C_1(\delta) = \rad(B)^{\delta-s}\cdot \max_{1\le j\le n(\delta)}
\left\{\frac{\prod_{i=0}^j R_i^\delta}{R_j^s}, 1\right\}
$$
is a constant independent of the choice of $E$. We continue with the
chain of upper inequalities to get
$$
\mu(E)\stackrel{\eqref{eq4}}\le C(X)\cdot C_1(\delta)\cdot C_2\cdot
(\rad(E))^{s-\delta},
$$
where
$$
C_2 = \max_{1\le j\le n_0}\left\{\prod_{i=0}^j
\frac{R_i^s}{t_i},1\right\}
$$
is again independent on the choice of $E$. By applying the Mass
Distribution Principle we conclude that $\dim\LKK(B,\RRR,\vs)\ge
s-\delta$. Since $\delta$ is arbitrary the lemma is proven. \endproof

\begin{lemma}[See Proposition 3
in~\cite{badziahin_velani_2011}]\label{lem_subcantor} Let
$\KKK({\rm I},\RRR,\vr)$   be as in Theorem~\ref{th_cantor2}. Then
there exists a local Cantor set  $$\LKK({\rm I},\RRR,\vs)
\subset \KKK({\rm I},\RRR,\vr), $$where
$$
\vs:=(s_n)_{n\in \ZZ_{\ge 0}}  \quad   { with  } \quad  s_n:=
\mbox{$\frac12$} \,  f(R_n) \, .
$$
\end{lemma}

\noindent{\it Proof of Theorem~\ref{th_cantor2}.}

By Lemma~\ref{lem_subcantor} we have that
$$
\dim \KKK(B,\RRR,\vr)\ge \dim  \LKK(B,\RRR,\vs) .
$$
Now fix some positive $s<\liminf\limits_{n\to\infty}(\dim
A_\infty(B)-\log_{R_n}\!2)$. Theorem~\ref{th3} gives us that for
every~$n$,
$$
\dim A_\infty(B) = \frac{\log f(R_n)}{\log R_n} =: d.
$$

 Then, there exists an integer $n_0$ such that
$$
s \, <  \, d-\log_{R_n}\!2  \quad {\rm \ for \  all }  \quad n>n_0
\, .
$$
Also  note that
$$
t_n=f(R_n)-s_n =\frac{f(R_n)}{2}
$$
and
$$
R_n^s<\frac{f(R_n)}{2}= t_n   \quad {\rm \ for \  all }  \quad
n>n_0 \, .
$$
Therefore, Lemma~\ref{lem_local} implies that
$$\dim \LKK({\rm I},\RRR,\vs)  \ge s   \, . $$
The fact that this inequality is true for any
$s<\liminf\limits_{n\to\infty}(d-\log_{R_n}\!2)$ completes the proof
of Theorem~\ref{th_cantor2}. \endproof

Finally we provide the theorem which shows that the intersection of
generalised Cantor sets on $X$ is often again a Cantor set.
\begin{theorem}[See Theorem~5 in~\cite{badziahin_velani_2011}]\label{th_icantor}
For each integer $1\le i \le k $, suppose we are given a generalised Cantor set $\KKK
(B,\RRR,\vr_i) $. Then
$$
\bigcap_{i=1}^{k} \KKK (B,\RRR,\vr_i)
$$
is a $(B,\RRR,\vr)$-Cantor set, where
$$
\vr:=(r_{m,n})  \quad\mbox{with } \quad r_{m,n} :=\sum_{i=1}^k
r^{(i)}_{m,n}  \, .
$$
\end{theorem}

With almost the same proof one can extend this theorem to countable
intersections of generalized Cantor sets.
\begin{theorem6}\label{th_icantorcountable}
For each integer $i\in\NN $, suppose we are given a generalised
Cantor set $\KKK (B,\RRR,\vr_i) $. Assume that the series
$$
r_{m,n}:= \sum_{i=1}^\infty r^{(i)}_{m,n}
$$
converges for all pairs $m,n\in \NN$ with $m\le n$ (or equivalently,
only finitely many of $r^{(i)}_{m,n}$ are non-zero). Then
$$
\bigcap_{i=1}^{\infty} \KKK (B,\RRR,\vr_i)
$$
is a $(B,\RRR,\vr)$-Cantor set with $ \vr:=(r_{m,n})$.
\end{theorem6}

\subsection{Images of generalized Cantor sets under bi-Lipschitz
map}

Let $\phi\;:\; X\to X$ be a bi-Lipschitz homeomorphism; i.e there
exists a constant $K>0$ such that
$$
\forall x_1,x_2\in X,\;\; K^{-1} \vd(x_1,x_2)\le \vd(\phi(x_1),
\phi(x_2))\le K\vd(x_1,x_2).
$$
One can easily check that then
$$
B(\phi(x), r/K)\subset \phi(B(x,r))\subset B(\phi(x), Kr).
$$
We denote the first (inscribed) ball by $I\phi(B)$ and the second (escribed) ball by $E\phi(B)$. We will also need a  slightly more restrictive packing condition than property (S4) enforced on the metric space $X$:
\begin{itemize}
\item[(S5)] For each $K\in\RR_{>1}$ there exists a constant $C(K,X)$
such that any ball $B$ of radius $Kr$ cannot intersect more than
$C(K,X)$ disjoint open balls of radius $r$.
\end{itemize}
One can easily check that property~(S4) of $X$ follows from
property~(S5) with $C(X) = C(3,X)$. Finally, note that the spaces $X$
appearing in examples (a) -- (c) from Section~\ref{sec1} satisfy condition~(S5).

\begin{theorem}\label{th8}
Let $(X,\SSS,U,f)$ be a splitting structure on a complete metric space $X$ satisfying condition~(S5). Assume also that  $A_\infty(B) = B$ for each ball
$B\in\BBB(X)$. Then for every bi-Lipschitz homeomorphism $\phi\;:\;
X\to X$ there exists a constant $C>0$ such that
$\phi(\KKK(B,\RRR,\vr))$ contains some $(I\phi(B),
\RRR,C\vr)$-Cantor set where
$$
C\vr:=\{Cr_{m,n}\;:\; m,n\in\ZZ_{\ge 0}, m\le n\}.
$$
\end{theorem}

\noindent {\bf Remark.} Surely the condition $A_\infty (B) = B$ is
quite restrictive. However, it is absolutely essential for the
theorem. One can check that the canonical splitting structures for
both $\RR^n$ and $\QQ_p^n$ satisfy that condition. On the other hand
the splitting structure $(\RR,\SSS,U,f)$ from example (c) does not
satisfy it.

\proof First, note that since $A_\infty (B) = B$ then for every
ball $B \in \BBB(X)$ and every $R\in U$ we have
$$
B = A_\infty (B)\subseteq
\bigcup_{b\in\SSS(B,R)} b\quad \Rightarrow\quad
\bigcup_{b\in\SSS(B,R)} b = B.
$$

Since $\KKK(B,\RRR,\vr)$ is a generalized Cantor set we have
collections $\BBB_n, \III_n$ and $\III_n^m$ (for $m< n$) associated with it
(see the Cantor set construction algorithm). We now outline the procedure for the construction of the generalised Cantor set inside $I\phi(B)$. Let
$\BBB_0^\phi:=\{I\phi(B)\}$. We next inductively construct a nested
collection $\BBB_1^\phi,\BBB_2^\phi,\ldots, \BBB_n^\phi,\ldots$.
Given a collection $\BBB_n^\phi$, construct the subsequent collection $\BBB_{n+1}^\phi$ via the following operations:
\begin{itemize}
\item Splitting procedure: Compute the collection
$$
\III_{n+1}^\phi:= \bigcup_{B_n^\phi\in\BBB_n^\phi} \SSS(B_n^\phi,R_n).
$$
\item Removing procedure: Remove all balls $B_{n+1}^\phi\in
\III_{n+1}^\phi$ for which
$$
\exists B_{n+1}\in \III_{n+1}\backslash \BBB_{n+1}\;\mbox{ s.t. }\;
B_{n+1}^\phi \cap \phi(B_{n+1}) \neq \emptyset.
$$
\end{itemize}
By construction we have that
$$
\bigcup_{B_n^\phi\in\BBB_n^\phi}B_n^\phi \subset
\phi\left(\bigcup_{B_n\in\BBB_n} B_n\right)
$$
and therefore the set
$$
\KKK^\phi:= \bigcap_{i=0}^\infty \bigcup_{B^\phi\in \BBB_i^\phi}
B^\phi
$$
is a subset of $\phi(\KKK(B,\RRR,\vr))$.

We will show that the set $\KKK^\phi$ is indeed $(I\phi(B), \RRR,
C\vr)$-Cantor for some constant $C>0$. Consider a ball $B_n^\phi \in
\BBB_n^\phi$. By construction, its radius is $K^{-1}\cdot
\prod_{i=0}^{n-1} R_i^{-1} \rad(B)$. If it intersects $\phi(B_n^*)$
for some $B_n^*\subset \BBB_n$ then it also intersects
$E\phi(B_n^*)$, whose radius is $$K\cdot \prod_{i=0}^{n-1} R_i^{-1}
\rad(B) = K^2\cdot \rad(B_n^\phi).$$ Therefore,
$$
\vd(\cent(B_n^\phi),\cent(E\phi(B_n^*)))\le (1+K^2)\rad(B_n^\phi).
$$
This in turn implies that
$$I\phi(B_n^*)\subset  (2+K^2)B_n^\phi.$$ Since
$$
\rad(I\phi(B_n^*)) = \rad(B_n^\phi),
$$
it follows from condition~(S5) that there are no more than
$C(K^2+2,X)$ balls $B_n^*\in \BBB_n$ such that $\phi(B_n^*)\cap
B_n^\phi \neq\emptyset$. By the same arguments we deduce that for a fixed ball $B_{n+1}\in
\BBB_{n+1}\backslash \III_{n+1}$ there are at most $C(K^2+2, X)$
balls $B_{n+1}^{\phi}\in \III_{n+1}^\phi$ which have nonempty
intersection with $\phi(B_{n+1})$.

Now we construct $\III^{n\phi}_{n+1}$ from $\III_{n+1}^\phi$ by
removing all balls $B_{n+1}^\phi\in \III_{n+1}^\phi$ which have
nonempty intersection with at least one of the sets $\phi(B_{n+1})$
where $B_{n+1}\in \III_{n+1}\backslash \III^n_{n+1}$. For a fixed
ball $B_n\in \BBB_n$ we have
$$
\# \{B_{n+1}^\phi\! \in \III_{n+1}^\phi: \exists B_{n+1} \in
\SSS(B_n, R_n)\cap (\III_{n+1}\backslash \III_{n+1}^n),
\phi(B_{n+1})\cap B_{n+1}^\phi \neq\emptyset\}\! \le C(K^2+2, X)
r_{n,n}.
$$
As we have already shown for a fixed $B_n^\phi \in \BBB_n^\phi$
there are at most $C(K^2+2,X)$ balls $B_n\in\BBB_n$ such that
$\phi(B_n)$ intersects $B_n^\phi$. Therefore, in total we have
$$
\# \{B_{n+1}^\phi \in \III_{n+1}^\phi \backslash \III_{n+1}^{n\phi}
\;:\; B_{n+1}^\phi \in \SSS(B_n^\phi,R_n)\} \le (C(K^2+2,X))^2
r_{n,n}.
$$

We proceed further with the Cantor set construction by constructing
the collection $I^{m\phi}_{n+1}$ from $I^{(m+1)\phi}_{n+1}$ ($0\le
m<n$) by removing all balls $B_{n+1}^\phi \in I^{(m+1)\phi}_{n+1}$
which have nonempty intersection with at least one of the sets
$\phi(B_{n+1})$ where $B_{n+1}\in \III_{n+1}^{m+1}\backslash
\III_{n+1}^m$. The same arguments as before yield for every ball $B_m^\phi \in\BBB_m^\phi$ the estimate

$$
\# \left\{B_{n+1}^\phi \in \III_{n+1}^{(m+1)\phi} \backslash
\III_{n+1}^{m\phi} \;:\; B_{n+1}^\phi \in
\SSS\left(B_m^\phi,\prod_{i=0}^{n-m}R_i\right)\right\} \le
(C(K^2+2,X))^2 r_{m,n}.
$$
This completes the proof that $\KKK$ is $(I\phi(B), \RRR,
C\vr)$-Cantor with $C:= (C(K^2+2,X))^2$.
\endproof

\section{Cantor-winning sets}\label{sec_winning}

Theorems~\ref{th_cantor1}~--~\ref{th_icantor} show that under
certain conditions on the sequences $\RRR$ and $\vr_i$ the finite
intersection
$$
\bigcap_{i=1}^{k} \KKK (B,\RRR,\vr_i)
$$
is non-empty or even has positive Hausdorff dimension. However they
do not cover countable intersections. Moreover, one can easily
provide a finite collection of generalized Cantor sets on $X$ which
have empty intersection. Theorem~6* on the other hand states that
under even stronger conditions we may deduce a similar statement for
a countable intersection of Cantor sets. However, all of these
conditions are somewhat cumbersome and may be quite difficult to
check. The aim of this section is to define a collection of sets
which satisfy properties (W1) and~(W2) of winning sets (that is, their
Hausdorff dimension equals to $\dim A_\infty(B)$ and their countable
intersection has the same property) and whose qualifying conditions
are much more clear cut.

Consider the constant sequence $\RRR = R,R,R,\ldots$. In this case
we will denote any associated generalised Cantor set $\KKK(B,\RRR,\vr)$ (respectively local Cantor set
$\LKK(B,\RRR,\vs)$) by $\KKK(B,R,\vr)$ (respectively
$\LKK(B,R,\vs)$). An easy inspection of the construction algorithm
for generalized Cantor sets on $X$ precipitates the following
proposition which will play a crucial role in constructing
our new class of sets.

\begin{proposition}\label{prop1}
Let $R\in U$, $k\in\NN$. Then $\KKK(B,R^k,\vr)$ is also $(B,R,\vt)$-Cantor where
\begin{equation}\label{eq_prop1}
t_{m,n} :=\left\{\begin{array}{ll} r_{m/k,(n+1)/k-1}&\mbox{if
}m\equiv n+1\equiv 0\pmod k;\\
0&\mbox{otherwise}.
\end{array}\right.
\end{equation}
\end{proposition}

Now, we are prepared to give the formal definition of Cantor-winning
set, the main object of interest in this paper.
\begin{definition}
Fix a ball $B\in\BBB(X)$. Given a parameter $\epsilon_0>0$ we say a
set $K\in X$ is \textbf{$\boldsymbol\epsilon_{\mathbf 0}$-Cantor-winning on $\mathbf B$ for the
splitting structure $\mathbf{(X,\SSS,U,f)}$} if for every
$0<\epsilon<\epsilon_0$ there exists $R_\epsilon\in U$ such that for every
$R\ge R_\epsilon$ with $R\in U$ the set $K$ contains some
$(B,R,\vr)$-Cantor set where
\begin{equation}\label{eq5}
r_{m,n} = f(R)^{(n-m+1)(1-\epsilon)}\quad\mbox{for every }m,n\in\NN,
m\le n.
\end{equation}
\end{definition}

If the splitting structure $(X,\SSS,U,f)$ is fixed then for
conciseness we omit its mention and simply say $K$ is
\textit{$\epsilon_0$-Cantor-winning on $B$}. Similarly, unless
otherwise specified a set $K\subset \RR^k$ or $K\subset \QQ_p^k$
will be referred to as being $\epsilon_0$-Cantor-winning on $B$ if
$K$ is $\epsilon_0$-Cantor-winning  on~$B$ with respect to the
relevant canonical splitting structure.

\begin{definition}
If a set $K\in X$ is $\epsilon_0$-Cantor-winning on~$B$ for every
ball $B\in\BBB(X)$ then we say that $K$ is \textbf{
$\boldsymbol\epsilon_{\mathbf 0}$-Cantor-winning}, and simply \textbf{Cantor-winning} if
$K$ is $\epsilon_0$-Cantor-winning for some $\epsilon_0>0$.
\end{definition}

We may apply Theorem~\ref{th_cantor2} to estimate the Hausdorff dimension
of Cantor-winning sets.
\begin{theorem}\label{th4}
If the complete metric space  $X$ satisfies condition~(S4). Then, for any $B \in \BBB(X)$ and any $\epsilon_0>0$ the Hausdorff dimension of an $\epsilon_0$-Cantor-winning set on $B$ is at least $\dim A_\infty(B)$.
\end{theorem}
\proof
If the splitting structure
$(X,\SSS,U,f)$ is trivial then $\dim \KKK(B,R,\vr) = \dim
A_\infty(B) = 0$. Otherwise, by taking if needed a power of $R$
in place of $R$ one can guarantee that $f(R)>4$. Also,
condition~\eqref{ineq_lem2} is obviously satisfied. Also, in this case the final condition~\eqref{cond_th2} condenses to the following:
$$
\sum_{k=0}^n f(R)^{(k+1)(1-\epsilon)} \left(\frac{4}{f(R)}\right)^k
\le \frac{f(R)}{4}.
$$
One can easily check that it is true for $f(R)$ large enough. So by again replacing $R$ with a proper integer power of $R$ if necessary we get that~\eqref{cond_th2} is satisfied. Thus, for any Cantor-winning set $E$ we have
$$
\dim E\ge \dim A_\infty(B) - \log_R 2.
$$
This estimate holds true with any integer power $R^k$ in place of $R$, and the theorem is proven.
\endproof
\begin{corollary}
    Let  $K$ be a Cantor-winning set. Then, for any $B \in \BBB(X)$ we have
    $$
    \dim(K \cap A_\infty(B)) \, = \, \dim(A_\infty(B)).
    $$
    In particular, in the case that $A_\infty(B)=B$ we have $\dim(K)=\dim(X)$.
\end{corollary}

 Next, we
 will show that the countable intersection of
 $\epsilon_0$-Cantor-winning sets is again $\epsilon_0$-Cantor-winning.

\begin{theorem}\label{th5}
Let a splitting structure $(X,\SSS,U,f)$ be nontrivial. Then, given
$\epsilon_0>0$ and a countable collection $\{K_i\}_{i\in\NN}$ of
$\epsilon_0$-Cantor-winning sets, the intersection
$$
\bigcap_{i=1}^\infty K_i
$$
is also $\epsilon_0$-Cantor-winning.
\end{theorem}

\proof

Consider an arbitrary positive $\epsilon<\epsilon_0$. By the
definition of $\epsilon_0$-Cantor-winning sets we have that $K_1$
contains $\KKK(B, R, \vr_1)$ for $R$ large enough where
$r_{m,n}^{(1)} = f(R)^{(n-m+1)(1-\epsilon)}$. Choose $R_\epsilon$ such
that $t\le f(R_\epsilon)^{t(\epsilon_0-\epsilon)}$ for any positive integer
$t$. Then for each $i>1$ one can inductively find $k_i\in\NN$ large
enough such that $k_{i+1}>k_i$ and the set $K_i$ contains
$\KKK(B,R_\epsilon^{k_i},\vr_i)$. Here $\vr_i$ are defined by the
formula~\eqref{eq5}:
$$
r_{m,n}^{(i)} = f(R_\epsilon)^{k_i(n-m+1)(1-\epsilon)}.
$$
By the definition of generalised Cantor sets any $(B,R_\epsilon,\vr_i)$-Cantor
set is also $(B,R_\epsilon,\tilde{\vr}_i)$-Cantor as soon as $r_{m,n}^{(i)}
\le \tilde{r}_{m,n}^{(i)}$. Therefore, without loss of generality we
can always assume that $\epsilon>\epsilon_0/2$. Next, we use
Proposition~\ref{prop1} to deduce that $K_i$ is also $(B,R_\epsilon,\vt_i)$-Cantor, where~$\vt_i$ is computed from $\vr_i$ by
formula~\eqref{eq_prop1}. This enables us to implement Theorem~$6^*$, which yields that
$$
\bigcap_{i=1}^\infty K_i \supset \bigcap_{i=1}^\infty
\KKK(B,R_\epsilon,\vt_i) = \KKK(B,R_\epsilon,\vt),
$$
where $t_{m,n}=\sum_{i=1}^\infty t_{m,n}^{(i)}$.

Finally, we must check that the values $t_{m,n}$ satisfy
condition~\eqref{eq5}. Notice that $t^{(1)}_{m,n}$ always contributes the value
$f(R_\epsilon)^{(n-m+1)(1-\epsilon)}$ to $t_{m,n}$. For $i>1$ this
contribution comprises
$$
f(R_\epsilon)^{k_i((n+1)/k_i - 1 - m)/k_i +1)(1-\epsilon)} =
f(R_\epsilon)^{(n-m+1)(1-\epsilon)}
$$
if $m\equiv n+1\equiv 0\pmod{k_i}$. Otherwise, $t^{(i)}_{m,n}$ does
not contribute anything to $t_{n,m}$. In other words, we have
\begin{equation}\label{eqn:setabove}
t_{m,n} = f(R_\epsilon)^{(n-m+1)(1-\epsilon)}\cdot \#\{i\in\NN\;:\; m\equiv
n+1\equiv 0\pmod{k_i}\}.
\end{equation}
Since all the numbers $k_i$ are distinct the cardinality of the set
on the right hand side of (\ref{eqn:setabove}) is at most $n-m+1$
and so we have
$$
t_{m,n}\le (n-m+1)f(R_\epsilon)^{(n-m+1)(1-\epsilon)}\le
f(R_\epsilon)^{(n-m+1)(1-2\epsilon+\epsilon_0)}.
$$
The final inequality holds due to the choice of $R_\epsilon$. As $\epsilon$
runs within the range $(\epsilon_0/2, \epsilon_0)$, the value
$2\epsilon-\epsilon_0$ takes any value within $(0,\epsilon_0)$.
Therefore, the intersection $\bigcap_{i=1}^\infty K_i$ contains a generalised Cantor set
$\KKK(B,R_\epsilon,\vt)$ satisfying property~\eqref{eq5}. Finally, the
same arguments apply if the parameter $R_\epsilon$ is replaced by any
other value of $R \in U$ with $R>R_\epsilon$. This completes the proof.
\endproof

\begin{theorem}
Let $(X,\SSS,U,f)$ be a splitting structure on a complete metric
space $X$ satisfying condition~(S5). Assume also that  $A_\infty(B)
= B$ for each ball $B\in\BBB(X)$ and let $\phi\;:\; X\to X$ be a
bi-Lipschitz homeomorphism. If $K\subset X$ is
$\epsilon_0$-Cantor-winning on a ball $B$ then $\phi(K)$ is
$\epsilon_0$-Cantor-winning on $I\phi(B)$.

\noindent Moreover if $K$ is $\epsilon_0$-Cantor-winning then so is
its image $\phi(K)$.
\end{theorem}

\proof If suffices to combine the definition of an
$\epsilon_0$-Cantor-winning set with Theorem~\ref{th8}. Indeed,
consider an arbitrary $0<\epsilon<\epsilon_0$. Then, by definition
there exists $R_\epsilon\in U$ such that for $R \geq R_\epsilon$ the set $K$ contains
$\KKK(B,R,\vr)$ where $r_{m,n}= f(R)^{(n-m+1)(1-\epsilon)}$. By
Theorem~\ref{th8}, the image $\phi(K)$ contains a $(I\phi(B),
R,C\vr)$-Cantor set for some absolute positive constant $C$
independent of $R$ and $\epsilon$. By choosing $\epsilon_1$
satisfying $\epsilon<\epsilon'<\epsilon_0$ and $R_{\epsilon'}$ large enough
so that $f(R_1)^{\epsilon'-\epsilon}>C$ it follows that for
$R>\max\{R_\epsilon,R_{\epsilon'}\}$ one has $Cr_{m,n}\le
f(R)^{(n-m+1)(1-\epsilon_1)}$. Thus, the set $\phi(B)$ is
$\epsilon_0$-Cantor-winning on $I\phi(B)$.

To prove the final statement we take an arbitrary ball $B\in
\BBB(X)$ and consider its preimage $\phi^{-1}(B)$. Take the escribed
ball $E\phi^{-1}(B)$. Since  $K$ is $\epsilon_0$-Cantor-winning it
is in particular  $\epsilon_0$-Cantor-winning on $E\phi^{-1}(B)$.
Therefore, the image $\phi(K)$ is $\epsilon_0$-Cantor-winning on
$I\phi(E\phi^{-1}(B))$. The final observation is that
$I\phi(E\phi^{-1}(B))=B$. This shows that $\phi(B)$ is indeed
$\epsilon_0$-Cantor-winning.
\endproof

{\bf Remark.} In~\cite{beresnevich_2013} the similar notion of
\textit{Cantor rich sets in $\RR$} was independently introduced. With
some effort this concept could also be generalised to $\RR^N$ and in turn arbitrary
complete metric spaces. Cantor rich sets are also known to satisfy
conditions~(W1) and~(W2). However, in the authors' opinion the conditions of
$\epsilon_0$-Cantor-winning sets are easier to check yet retain the same desirable properties. Furthermore, the following section provides some reasoning as to why our setup may be preferable in many cases (see the Remark at the close of \S\ref{McM} and~\cite{BHN}). It would
be interesting to compare the two notions, to ask whether the two
concepts are equivalent, whether one of them includes another, or if
neither of these two possibilities hold, although this appears to be a quite difficult and nuanced question.

\section{Relationship with classical winning sets}\label{sec:relationship}

We have shown that under certain conditions Cantor-winning sets
satisfy the same desirable properties (W1) -- (W3) as classical
winning sets in $\RR^N$. It is therefore natural to ask if and how
these two concepts are compatible.

%
%
%

In his original paper, Schmidt defined his game in the context of
any complete metric space~$X$. For the $(\alpha, \beta)$-game played
on $X$, Alice and Bob pick successive nested balls in the same
manner as described in Section \ref{sec:winning}. The definitions of
$\alpha$-winning sets and winning sets for gameplay in an arbitrary
complete metric space are entirely analogous to those for the game
played in $\RR^N$. Strictly speaking, since a generic ball in $X$
may not necessarily have a unique  centre or radius, Alice and Bob
should pick successive pairs of centres and radii satisfying some
partial ordering
as opposed to simply picking successive nested balls. However, for
the sake of clarity one may simply assume that this nuance is
accounted for in each of Alice's and Bob's strategies.

Properties (W2) and (W3) are satisfied by winning sets for any
$(\alpha, \beta)$-game played in an arbitrary complete metric space
$X$ (see \cite{schmidt_1966} and \cite{dani_1989} respectively). On
the other hand,  winning sets need not satisfy  property (W1) in
general. Indeed, Proposition $5.2$ of \cite{kleinbock_weiss_2010}
provides an example of a winning set of zero Hausdorff dimension. In
\cite{kleinbock_weiss_2010} it is also shown that if $X$ supports a
measure satisfying certain desirable rigidity properties then
property (W1) does indeed hold. We discuss one such property in a
later section - see (\ref{powerlaw}).

Comparing directly the property of being a winning set in $X$ with
the property of  being a Cantor-winning in $X$ appears to be a very
difficult problem and would likely require lengthy and technical
discussion. For this reason, and to help maintain the flow of this
paper, we only mention that the authors intend to return to this
topic in the subsequent work \cite{BHN}.
It is though much more feasible to place our framework within the hierarchy of various classes of games related to those of Schmidt that exhibit a
slightly higher level of rigidity. Indeed, as we will see the
relationship between Cantor-winning sets and the `winning sets' of
these classes of games is rather clear cut.

\subsection{McMullen's Game} \label{McM}
In \cite{mcmullen_2010}, McMullen proposed the following
one-parameter variant of Schmidt's game, defined in such a way that
instead of choosing a region where Bob must play, Alice must now
choose a region where he must not play. To be precise, first choose
some parameter $\beta \in (0, \gamma(X))$, where $\gamma(X)>0$ is
some absolute constant  (to be determined later) depending on the
metric space $X$. McMullen's \textit{$\beta$-absolute game} begins
with Bob picking some initial ball $B_1 \in \BBB(X)$. Alice and Bob
then take it in turns to place successive balls in such a way that
$A_i \subset B_i$ and
$$
B_1 \: \supset\: B_1 \setminus A_1 \:\supset \:B_2 \:\supset\: B_2 \setminus A_2 \:\supset \:B_3 \:\supset \: \cdots,
$$
subject to the conditions
$$
\rad(B_{i+1}) \geq \beta \cdot \rad(B_i), \quad \rad(A_{i+1}) \leq \beta \cdot \rad(B_i), \quad \forall \, i \in \NN.
$$
We say a set $E \subset X$ is
\textit{$\beta$-absolute winning} if Alice has a strategy which guarantees
\begin{equation}\label{absolutewinning}
\bigcap_{i \in \NN} \, B_i \cap E \, \neq \, \emptyset
\end{equation}
for the game with parameter $\beta$. The set $E$ is said to be
\textit{absolute winning} if it is $\beta$-absolute winning for
every $\beta\in(0,\gamma(X))$. Note that in general $\bigcap_i B_i$ may not
necessarily be a single point as in Schmidt's $(\alpha,
\beta)$-game.

McMullen's original definition of the $\beta$-absolute game
exclusively involved the  selection of closed balls in $\RR^N$.
However, the mechanics  described above make sense when outlining
the rules for play with (metric) balls in any complete metric
space~$X$.

The purpose of the upper bound $\gamma(X)$ for the choice of
$\beta$, as introduced above, is to ensure that at every stage of
a $\beta$-absolute game there is always a legal place for Bob to
place his ball wherever Alice may have placed her preceding ball.
For the game played on $X=\RR^N$ with Euclidean balls one may take
$\gamma(X)= 1/3$ as per McMullen's original definition.  To see that
this condition is necessary, notice that for $\beta \geq 1/3$ Alice
may then at any stage choose her ball $A_i$ to simply be the ball
$B_i$ scaled down by $\beta$. In doing so she would leave no
possible choice of ball $B_{i+1}$ satisfying  $B_i\setminus A_i
\:\supset \:B_{i+1}$. However, for $\beta<1/3$ such a choice is
always possible in $\RR^N$.

For the game played on an arbitrary complete metric space such a constant $\gamma(X)$ need not exist.
However, it was recently observed in
\cite{fishman_simmons_urbanski_2013} (see their Lemma $4.2$, and
also \cite{mayeda_merrill_2013, Weil}) that a sufficient condition
for the existence of $\gamma(X)$ is  that the metric space in
question is \textit{uniformly perfect}. Recall that for $0<c<1$ a
metric space $X$ is said to be\textit{ $c$-uniformly perfect} if for
every metric ball $B(x,r) \neq X$ we have $B(x, r) \setminus  B(x,
cr) \neq \emptyset$, and is said to be \textit{uniformly perfect} if
it is $c$-uniformly perfect for some $c$. If the metric space $X$ is
uniformly perfect one may then take $ \gamma(X) = c/5$, although it
should be noted that this is not necessarily the optimal (largest
possible) choice. It is easy to see that if $X$ is endowed with a
non-trivial splitting structure and satisfies condition $(S4)$ with
constant $C(X)$ then $X$ is indeed uniformly perfect and so in the
setting of this paper McMullen's game is always well defined. In
particular, one may take $c=u_0^{-1}$, where $u_0 \in U$ is the
smallest natural number for which $f(u_0)>C(X)$.

It is well known that an absolute winning set in $\RR^N$ is
$\alpha$-winning for every $\alpha \in (0,1/2)$, and that for the
game played on a $c$-uniformly perfect complete metric space an
absolute winning set is $\alpha$-winning for every $\alpha \in
(0,c/5]$ - see \cite{mcmullen_2010} and
\cite{fishman_simmons_urbanski_2013} respectively. In both cases it
can be shown that the countable intersection of $\beta$-absolute
winning sets is again $\beta$-absolute winning, and that the image
of an absolute winning set under a bi-Lipschitz homeomorphism is
again absolute winning. Thus, absolute winning sets also satisfy
properties (W2) and~(W3). In fact, it is the case (see Proposition
$4.3(v)$ of \cite{fishman_simmons_urbanski_2013}) that absolute
winning sets in any uniformly perfect complete metric space satisfy
the following slightly stronger version of the latter property:
\begin{itemize}
    \item[(W$3^\ast$)] The image of any absolute winning set under a quasisymmetric  homeomorphism  is
    again absolute winning.
\end{itemize}
As before, absolute winning sets do not in general satisfy condition
(W1), although if $X$ supports a measure satisfying (\ref{powerlaw})
it follows that property (W1) does hold. See
\cite{fishman_simmons_urbanski_2013} and the references therein for
further criterion.

The following theorem reveals that absolute winning sets have an
extremely clear cut relationship with Cantor-winning sets. We delay
the proof to a later subsection.

\begin{theorem}\label{thmabsolutewinning}
Assume a complete metric space $X$ is endowed with a non-trivial
splitting structure~$(X,\SSS,U,f)$ and that condition (S4) holds
with constant $C(X)$. If $E \subset X$ is absolute winning then $E$
is $1$-Cantor-winning.
\end{theorem}

{\bf Remark.} Since completion of this project, the authors (in collaboration with with Nesharim) \cite{BHN} have been able to show that the converse statement is in fact true, at least in the case of $\RR^N$ with canonical splitting structure. That is; remarkably, the property of being $1$-Cantor-winning in $\RR^N$ is in fact equivalent to the property of being absolutely winning!

\subsection{The Hyperplane Absolutely Winning game and its variants}
In \cite{bfkrw_2012}, a class of variants of McMullen's game was
introduced, the so-called \textit{$k$-dimensional absolute winning}
games. The most commonly utilised of these games is the
\textit{hyperplane absolute winning} (or \textit{HAW}) game. The
class of games in \cite{bfkrw_2012} was specifically defined for
play on subsets of $\RR^N$ and relies upon the existence of an
underlying vector space. For this reason, in order to discuss
$k$-dimensional absolute winning games in the full setting of this
paper we would first have to attach further structure to our
complete metric space $X$. In particular, if so inclined one could
define the games for subsets of some given Banach space, but since
such an extension has not yet appeared in the literature we content
ourselves with the setting of $\RR^N$ (with metric
$\mathbf{d}(x , y) = |x-y|_\infty$ and canonical splitting
structure) for the sake of clarity. Accordingly, we will refer to
the metric balls in $\RR^N$ as `boxes'. That said, one should
observe that an analogues method to the one we shall exhibit would
be applicable to questions concerning $k$-dimensional absolute
winning games played on more exotic spaces.

Firstly, fix $k \in \left\{0, 1,\ldots , N-1\right\}$ and some
parameter $0 < \beta < 1/3$. The $k$-dimensional  $\beta$-absolute
winning game has the same premise as McMullen's game in that Alice
must choose a region where Bob must not play, only now that region
is defined by the neighbourhood of a $k$-dimensional hyperplane of
$\RR^N$ rather than the neighbourhood of a single point. The game
begins with Bob picking some initial box~$B_1 \subset \RR^N$. Now,
assume Bob has played his $i$-th box~$B_i$. Then the game proceeds
with Alice choosing $\delta_i\le \beta$ and an affine subspace
$\mathcal L_i$ of dimension $k$ and removing its $(\delta_i \cdot
\rad(B_i))$-neighbourhood
$$
\mathcal{L}_i^{(\delta_i \cdot \rad(B_i) )}=\left\{x \in \RR^N: \:
\inf_{y \, \in \, \mathcal{L}_i}\mathbf |x-y|_\infty < \delta_i\cdot
\rad(B_i) \right\}
$$
from the box $B_i$. In accordance with this procedure, set $A_i:=
\mathcal{L}_i^{(\delta_i \cdot \rad(B_i))} \cap B_i$. Then, Bob for
his $(i+1)$-th move may choose any box $B_{i+1} \subset B_i
\setminus A_i$ satisfying $\rad(B_{i+1}) \geq \beta \cdot
\rad(B_i)$. A set $E \subset \RR^N$ is said to be
\textit{$k$-dimensionally $\beta$-absolute winning} Alice has a strategy guaranteeing that
$$
\bigcap_{i=1}^\infty \, B_i \: \cap \: E \: \neq \: \emptyset
$$
for the game played with parameter $\beta$. We simply say that $E$
is \textit{$k$-dimensionally absolute winning} if it is
\textit{$k$-dimensionally $\beta$-absolute winning} for every
$\beta \in (0, 1/3]$.

In the weakest case `$k=N-1$', the game is often referred to as the
\textit{hyperplane absolute winning} game for obvious reasons. For
simplicity, an $(N-1)$-dimensionally absolute winning set is then
referred to as being \textit{hyperplane absolute winning (HAW)}. One
can readily observe that the strongest case `$k=0$' coincides with
McMullen's original game on $\RR^N$. For every $k$, if a set is
$k$-dimensionally absolute winning sets then it is $\alpha$-winning
with respect to Schmidt's game for any $\alpha \in (0, 1/2)$. We
direct the reader to \cite{bfkrw_2012} for further discussion of the
properties of $k$-dimensionally absolute winning sets.

\begin{theorem}\label{thmhyperplanewinning}
Assume a subset $E \subset \RR^N$ is  $k$-dimensionally absolute
winning for some integer $k \in \left\{0, 1,\ldots , N-1\right\}$.
Then, the set $E$ is  $ \frac{N-k}{N}$-Cantor-winning.
\end{theorem}
Note that the case `$k=0$' corresponding to McMullen's game is
contained within the statement of Theorem~\ref{thmabsolutewinning}.
Broadly speaking, Theorems \ref{thmabsolutewinning} \&
\ref{thmhyperplanewinning} demonstrate that the property of being a
Cantor-winning set is weaker than the property of being a
$k$-dimensionally absolute winning set (for any given $k \in
\left\{0, 1,\ldots , N-1\right\}$). However, this weakening does not
come at the cost of losing properties (W1)~-~(W3).

\subsection{Proof of Theorems \ref{thmabsolutewinning} \& \ref{thmhyperplanewinning}  }

\subsubsection{Preliminaries}
In order to present our proofs we first require some terminology. For consistency we use the notation originally introduced in~\cite{schmidt_1966}. Additionally, for $k = 0, 1, \ldots, N-1$ let $\mathcal H_k$ denote the set of all affine $k$-dimensional hyperplanes in $\RR^N$.

In
each of the $k$-dimensionally absolute winning games (including
McMullen's game on a metric space $X$),  a set $E$ is
($k$-dimensionally) absolute winning if however we choose to place
Bob's  balls $B_i$ Alice has a `strategy' for placing her moves $A_i$
so that the set $\bigcap_i B_i$ intersects~$E$. Formally, a \textit{strategy} $F:=(f_1, f_2, \ldots)$ is a sequence of functions $f_i: \BBB(X)^i \rightarrow \mathcal H_k \times \RR_{>0}$.
Given a fixed parameter $\beta$, we say a strategy $F$ is \textit{legal} for the ($k$-dimensional) $\beta$-absolute game if it satisfies the following property for any finite sequence $(b_1, \ldots, b_n)$ of balls, any hyperplane~$h$, and any $s \in \RR_{>0}$:
 \begin{equation}\label{legal}
 \mbox{if } f_n(b_1, \ldots, b_n)=(h,s), \mbox{ then }
 s \leq \beta \cdot \rad(b_n).
 \end{equation}
 For $(h, s) \in \mathcal H_k \times \RR_{>0}$ denote by $g(h,s):=h^{(s)}$ the standard closed $s$-neighbourhood of the hyperplane $h$. We say $F$ is a \textit{winning strategy (for $E$)} with respect to the game with parameter $\beta$ if firstly it is legal and secondly if it then determines where Alice should place her moves $A_i:=g(f_i(B_1, B_2, \ldots, B_i))$ in such a way that, however we choose to place Bob's balls $B_1$ and $B_{i+1} \subset B_{i}\setminus A_{i}$ (for $i \in \NN$) in the game, condition~(\ref{absolutewinning}) holds. It is easily verified that a set $E$ is ($k$-dimensionally) $\beta$-absolute winning if and only if there exists a winning strategy for $E$ with respect to the ($k$-dimensionally) $\beta$-absolute game. 
 

The following key observation made by Schmidt in \cite{schmidt_1966} (his Theorem $7$) allows us to significantly simplify our notation: In any of the above games, the existence of a winning strategy for a set $E$ guarantees the existence of a `positional' winning strategy for $E$. We say a winning strategy~$F:=(f_1, f_2, \ldots )$ is \textit{positional} if each function $f_i$ depends  only upon the ball in its final component; that is, the placement of each of Alice's moves in the winning strategy depends only upon the position of Bob's immediately preceding ball, not on the entirety of
the game so far. For this reason, if the ball $b$ appears as Bob's $n$-th move during gameplay then we will without loss of generality write $g(f_n(b))$ to denote Alice's subsequent move as determined by the strategy~$F$.

As a final piece of terminology from~\cite{schmidt_1966}, given a target set $E$ we refer to a sequence~$(B_1, B_2, \ldots)$ of balls as an \textit{$F$-chain} if it consists of the moves Bob has made during a ($k$-dimensional) $\beta$-absolute game in which Alice has followed the winning strategy $F$ for~$E$. By definition we must have  (\ref{absolutewinning}) holds for this sequence. Furthermore, we say a finite sequence~$(B_1, B_2, \ldots B_n)$ is an \textit{$F_n$-chain} if there exist~$B_{n+1}, B_{n+2}, \ldots$ for which the infinite sequence~$(B_1, B_2, \ldots B_n, B_{n+1} \ldots )$ is an $F$-chain.

\subsubsection{Proof of Theorem \ref{thmabsolutewinning}}

Recall that for any non-trivial splitting structure satisfying condition ($S4$) the quantity $u_0\in U$ is defined to be the
smallest number such that $f(u_0)>C(X)$. By assumption our set $E\subset X$ is $\beta$-absolute winning for every $\beta < \gamma(X):=(5u_0)^{-1}$. 
Fix some ball $B \in \BBB(X)$ and $\epsilon \in (0,1)$, and let $R_1 \in U$ be the smallest integer for which $5u_0 < R_1$. Next, choose $R_2 \in U$ large enough so that for any $R\in U$ with $R \geq R_2$ we have $ f(R)^{(1-\epsilon)} \geq C(X)$. This is always possible for a non-trivial splitting structure by Corollary~\ref{corol1} and the multiplicativity of~$f$. Now set $R_\epsilon: = \max(R_1, R_2)$. 
To prove the theorem  it suffices to construct for each  $R\in U$ with $R \geq R_\epsilon$ a local Cantor set~$\mathcal{LK}(B, R, \vs)$ lying
inside $E$  for which $s_{n}\leq
f(R)^{(1-\epsilon)}$.

Fix some $R \in U$ satisfying $R \geq R_\epsilon$. Our method for constructing the set $\mathcal{LK}(B, R, \vs)$ is as follows. We play as Bob in an iteration of McMullen's game
with parameter $\beta= 1/R$. By assumption the set $E$ is 
$(1/R)$-absolute winning and so there exists a winning strategy~$F$ for~$E$ with respect to this game. Here, we have $\mathcal{H}_0= X$ and so $\{ g(h,s): \, (h, s) \in X \times \RR_{>0} \}$ coincides with the set of all closed balls $\BBB(X)$ in $X$.

Assume that Bob plays his first ball in position $B_1=B$ and allow the strategy $F=(f_1, f_2, \ldots)$ to determine Alice's first ball~$A_1:=g(f_1(B))$. Consider the set~$\SSS(B,R)$. Since by~(\ref{legal}) we have $\rad(A_1) \leq \frac 1 R \rad(B)=\rad(b)$ for every $b \in \SSS(B,R)$, the ball
$A_1$ may intersect at most $C(X)$ balls from the collection
$\SSS(B,R)$. 

The construction of the local Cantor set $\mathcal{LK}(B, R, \vs)$ comprises the construction of subcollections $\BBB_i
\subset \SSS(B, R^i)$ and a sequence $\vs=(s_n)_{n\in \ZZ_{\ge 0}}$. As a first step in this
procedure, define $\BBB_0:= \{ B \}$ and
$$
\BBB_1: =  \left\{ b \in \SSS(B,R): \: g(f_1(B)) \cap b = \emptyset \right\}.
$$
Upon setting $s_{0}:=\#(\SSS(B,R) \setminus \BBB_1  )$ we have $s_{0} \leq C(X) \leq f(R)^{(1-\epsilon)}$ as
required. Furthermore, any ball $B_2 \in \BBB_1$ is a legal choice for Bob's next move in the game; i.e., the finite sequence~$(B, B_2)$ is an $F_2$-chain for any $B_2 \in \BBB_1$. 

Assume now that for some $n \in \NN$ we have constructed the collections $\BBB_i$ and defined the values $s_{i-1} $  for $i=1, \ldots, n$. Assume also that these collections satisfy the property that for every $b \in \BBB_i$ we have $g(f_i(b')) \cap b = \emptyset$, where $b'$ is the unique ball in the collection $\BBB_{i-1}$ containing $b$. It is immediate that any finite sequence $(B_1, \ldots B_{n+1} )$ with $B_i \in \BBB_{i-1}$ is an $F_{n+1}$-chain. We construct the collection~$\BBB_{n+1}$ in the following way. Simply notice that for any $b' \in \BBB_n$ the ball~$g(f_{n+1}(b')$ may, by~(\ref{legal}) and condition~($S4$), intersect at most $C(X)$ of the balls from the collection $S(b', R)$. Indeed, for $b' \in \BBB_{n}$ let
$$
\BBB^{b'}_{n+1}: =  \left\{ b \in \SSS(b',R): \:  \: g(f_{n+1}(b')) \cap b = \emptyset \right\},
$$
and set
$$
\BBB_{n+1}: = \bigcup_{b' \in \BBB_n} \, \BBB^{b'}_{n+1} \quad\quad
\text{and } \quad\quad s_{n}:= \displaystyle\max_{b' \in
\BBB_n}\#\left(\SSS(b',R) \setminus \BBB^{b'}_{n+1}\right).
$$
Then, it follows that $ s_{n} \leq C(X) \leq f(R)^{(1-\epsilon)}$ and by definition that for every $b \in \BBB_{n+1}$ we have $g(f_{n+1}(b')) \cap b = \emptyset$, where $b'$ is the unique ball in the collection $\BBB_{n}$ containing $b$. Furthermore, if $(B_1, \ldots B_{n+1} )$ with $B_i \in \BBB_{i-1}$ is an $F_{n+1}$-chain then $(B_1, \ldots B_{n+1}, b )$ is an $F_{n+2}$-chain for any $b \in \BBB_{n+1}$.  This completes the inductive procedure. 

Upon defining
$$
\mathcal{LK}(B, R, \vs): \, = \, \bigcap_{i \in \ZZ_{\geq0}} \: \bigcup_{b \in \BBB_i} b,
$$
it only remains to show that $ \mathcal{LK}(B, R, \vs) \subseteq E$. With this in mind, choose some point $x \in \mathcal{LK}(B, R, \vs)$ and let $s=(b_i)_{i \in \NN}$ with $b_i \in \BBB_{i-1}$ be a sequence of balls for which $\bigcap_{i \in \NN}b_i = x$. By construction, we have ensured that each finite subsequence $(b_1, \ldots, b_n )$ is an $F_n$-chain. Moreover, it is readily verified that if $(b_1, b_2, \ldots)$ is a sequence of balls such that for every $n \in \NN$ the finite sequence $(b_1, \ldots b_n)$ is an $F_n$-chain, then $(b_1, b_2, \ldots )$ is an $F$-chain (c.f. \cite[Lemma $1$]{schmidt_1966}). It follows that condition~(\ref{absolutewinning}) holds and, since $\rad(b_i) \rightarrow 0$ as $i \rightarrow \infty$ implies the intersection $x=\bigcap_{i \in \NN}b_i$ is a single point, that $x \in E$ as required.

Finally, by the fact that the initial ball $B$ and the quantity $\epsilon \in (0,1)$ were arbitrary it follows that the set~$E$ is $1$-Cantor-winning.

\subsubsection{Proof of Theorem \ref{thmhyperplanewinning}}

The proof follows very similar arguments to those in the proof of Theorem~\ref{thmabsolutewinning}. For this reason we
only outline the modifications required. The key observation is that, given any box $B$
of sidelength $\diam(B)$ in $\RR^N$ and any $R \in \NN$,  the
rectangular neighbourhood
$$
\mathcal{L}^{(\rad(B)/R )} = \left\{x \in \RR^N: \:
\inf_{y \, \in \, \mathcal{L}}\mathbf |x-y|_\infty \, < \, \frac {\rad(B)} R \right\}
$$
of any $k$-dimensional hyperplane $\mathcal{L}$ passing through $B$
may intersect at most $c(k, N) \cdot R^{k}$ of the boxes $b \in S(B,
R)$. Here, the quantity $c(k, N) \in \RR_{>0}$ is an absolute
constant depending only upon $k $ and $N$.

Suppose the set $E\subset \RR^N$ is $k$-dimensionally absolute winning. Fix some box $B \subset \RR^N$ and some $\epsilon \in (0,(N-k)/N)$. Next, choose $R_\epsilon > 3$ large enough so that for any $R \geq R_\epsilon$ we have $ R^{N(1-\epsilon)} \geq c(k, N) \cdot R^k$. This is always possibly since $$1-\epsilon \: \:> \: \:1-\frac {N-k}{N} \:\:=\:\: \frac k N.$$  

Fix some $R \geq R_\epsilon$ and set consider a $\beta$-absolute game with $\beta= 1/R$. By assumption there exists a winning strategy $F=(f_1, f_2, \ldots,)$ associated with the parameter $\beta$ and the set $E$. Let Bob initially play the box $B_1:=B$ and set $\BBB_0:=\{ B \}$. 
As in the proof of Theorem \ref{thmabsolutewinning} one must
construct collections $\BBB_n \in S(B, R^n)$ and a sequence $\vs=(s_i)_{i \in \NN}$ in an iterative
fashion in order to define some local Cantor set $\mathcal{LK}(B, R, \vs)$.  Given $n \geq 0$, for every box $b\in\BBB_{n}$ played by Bob during gameplay the strategy $F$ determines the position and a neighbourhood of a $k$-dimensional affine hyperplane instructing Alice where to play her next move. By the above observation any such neighbourhood may intersect at most
$c(k, N)\cdot R^{k}$ boxes from the collection $S(b, R)$. Following
exactly the method of Theorem \ref{thmabsolutewinning}
one may analogously construct the required collections
$$
\BBB_{n+1}^{b'}:= \{b\in\SSS(b',R)\;:\; g(f_{n+1}(b'))\cap
b=\emptyset\}\quad\mbox{ for } b' \in \BBB_n \quad \mbox{ and }\quad \BBB_{n+1}:=\bigcup_{b'\in
\BBB_{n}} \BBB_n^{b'}.
$$
Furthermore, we may choose $s_{n}:= c(k, N)\cdot R^{k} \leq R^{N(1-\epsilon)}$ as required.
As before, since $F$ is a winning strategy it is ensured that the resulting local $(B, R, \vs)$-Cantor set falls inside~$E$. 

\section{Generalized badly approximable sets}\label{sec_genbad}

In~\cite{kristensen_thorn_velani_2006} the authors introduced a
broad notion of badly approximable sets. We now discuss how their
setup is related to ours, and more importantly, how we are able to
generalise their results. We begin by giving a brief outline of the
framework outlined in~\cite{kristensen_thorn_velani_2006}, tailored
to our needs.

Let $X$ be a complete metric space (with metric $\vd$) and let $\cR$
be a family of subsets $\cR:=\{R_\alpha\subset X\;:\; \alpha\in S\}$
indexed by an infinite countable set $S$. In  most of the
applications we will discuss, the subsets $R_\alpha$ will consist
simply of points in $X$. The subsets  $R_\alpha$ will be referred to
as {\it resonant sets}. We attach a `weight' to each resonant set by
introducing a function $h\;:\; S\to \RR_{\ge 0}$. For convenience we
will always assume that $h$ is bounded above by some absolute
constant; in other words, there exists a constant $C>0$ such that
for every $\alpha\in S$ we have $h(\alpha) \le C$. Next, for any set
$R\subset X$, let
$$
\Delta(R,\delta):= \{\vx\in X\;:\; \vd(\vx,R)\le \delta\}
$$
denote the  $\delta$-neighbourhood of $R$. Finally, we say a set of the form
$$
\bad(\cR,h):= \{\vx\in X\;:\; \exists c>0, \forall \alpha\in S,
x\not\in \Delta(R_\alpha, c\cdot h(\alpha))\}.
$$
is a \textit{generalised bad set}.

{\bf Remark.} Our definition slightly differs to that given
in~\cite{kristensen_thorn_velani_2006}. For simplicity we have
combined the two functions $\beta_\alpha$ and $\rho$ present
in~\cite{kristensen_thorn_velani_2006} into one function $h$. With
reference to the notation of~\cite{kristensen_thorn_velani_2006}, if
we take $\beta_\alpha = (h(\alpha))^{-1}$, $\rho(x) = x^{-1}$ and
$\Omega = X$ then $\bad^*(\cR,\beta,\rho)$ defined
in~\cite{kristensen_thorn_velani_2006} is precisely our set
$\bad(\cR,h)$. Furthermore, since there is a bijection between $S$
and our sequences $\cR$, in applications we will often use the
notation $h(R)$ for $R\in\cR$ instead of $h(\alpha)$.

The following all provide basic examples of generalised bad sets:
\begin{itemize}
\item The standard set $\bad$ of badly approximable numbers.
In this case an easy inspection shows that $\bad = \bad(\cR,h)$
where $\cR$ consists of rational points and $h(p/q):= 1/q^2$ for
every $p\in\ZZ, q\in\NN, \gcd(p,q)=1$.

\item The set $\bad_N$ of badly approximable points in $\RR^N$ (see subsection \ref{sec:winning} for precise definition).
Again one can check that $\bad_N$ is a generalised bad set for $$\cR =
\{\vp/q\;:\; \vp\in\ZZ^N, q\in\NN, \gcd(p_1,\ldots,p_N,q)=1\}$$ and
$h(\vp/q) = q^{-1-1/N}$.
\item The set $\bad_p:=\bad_p(1)$ of \textit{$p$-adically badly approximable numbers} as defined in~\eqref{def_badp}.
The implicit inequality in the definition of $\bad_p$ is clearly
satisfied for $q=0$. We can also without loss of generality assume
that $\gcd(q,r)=1$. Then, by dividing both sides of the inequality
in~\eqref{def_badp} by $|q|_p$ one can check that it is a
generalised bad set for $\cR = \QQ\subset \QQ_p$ and $$h(r/q) =
(\max\{|r|^2,|q|^2\}\cdot |q|_p)^{-1}.$$
\end{itemize}

In~\cite{kristensen_thorn_velani_2006} the authors give quite
general conditions on $\cR$ and $h$ which guarantee that a generalised bad  set $\bad(\cR,
h)$ has full Hausdorff dimension. Namely they prove the following.

\begin{theoremKTV}[Theorem 1 in~\cite{kristensen_thorn_velani_2006}]
Let $X$ support a measure $m$ for which there exist strictly
positive constants $\delta$ and $r_0$ such that for any $x\in X$ and
$r\le r_0$,
\begin{equation}\label{powerlaw}
ar^\delta\le m(B(x,r))\le br^\delta,
\end{equation}
where $0<a\le 1\le b$ are constants independent of the ball. Define
$J(n):= \{\alpha\in J\;:\; R^{n-1}\le (h(\alpha))^{-1}<R^n\}$.

Assume that for $R$ large enough there exists $\theta\in\RR^+$ so
that for $n\ge 1$ and any ball $B_n = B(x,h(R^n))$ there exists a
collection $C(\theta B_n)$ of disjoint balls such that
$$
\forall B_{n+1}\in C(\theta B_n),\; \rad (B_{n+1}) = 2\theta
h(R^{n+1});\; B_{n+1}\subset B(x,\theta h(R^n));
$$$$
\#C(\theta B_n) \ge \kappa_1 R^\delta
$$and$$
\#\{B_{n+1}\in C(\theta B_n)\;:\; \exists \alpha\in J(n+1)\mbox{
s.t. } \cent(B_{n+1})\in \Delta(R_\alpha, 2\theta h(R^{n+1}))\}\le
\kappa_2 R^\delta,
$$
where $0<\kappa_2<\kappa_1$ are absolute constants independent of
$k$ and $n$. Furthermore, suppose that $\dim (\cup_{\alpha\in J}
R_\alpha)<\delta$. Then, $\dim (\cR,h)=\delta = \dim X$.
\end{theoremKTV}

This theorem provides the Hausdorff dimension for $\bad(\cR, h)$ in
a wide ranging setup in which relatively mild (but rather technical)
conditions on $\cR, X$ and $m$ are assumed. However, we show that
some sets $\bad(\cR, h)$ which do not satisfy certain conditions of
Theorem~KTV can still be shown to satisfy properties (W1)~--~(W3),
as do many sets which do fall within the scope
of~\cite{kristensen_thorn_velani_2006}. In this sense our framework
is more far reaching than that presented
in~\cite{kristensen_thorn_velani_2006}. On the other hand, in order
to do this we will need to impose slightly more structure on the
balls $B_n$ and classes $C(\theta B_n)$, which in turn makes some of
our conditions slightly stronger than those imposed in Theorem~KTV.

%

Observe that one may consider the sets $\bad(\cR,h)$ as the set of
points surviving after the removal of every neighborhood
$\Delta(R_\alpha,c\cdot h(\alpha))$ from $X$. On adopting this point
of view one may appreciate the similarity between general bad sets
and generalized Cantor sets. To further illustrate this connection
we now provide an algorithm, which will be referred as a {\it bad to
Cantor set construction}, demonstrating that the intersection of
every generalised bad set $\bad(\cR,h)$  with any set $A_\infty(B)$
contains some generalized Cantor set $\KKK(B,R,\vr)$.

{\bf Bad to Cantor Set Construction:}
\begin{enumerate}
\item Fix $R$ large enough and choose $c$ small enough such that
\begin{equation}\label{eq_rc}
\sup_\alpha\{c\cdot h(\alpha)\} \le \diam(B)\cdot R^{-1}.
\end{equation}
This can be done since $h(\alpha)$ is always bounded above by an
absolute constant.

\item Split the collection $\cR$ into classes $C(n)$, for $n\in\NN$, in
the following way. Let
\begin{equation}\label{eq_cn}
C(n):=\{R_\alpha\in\cR\;:\; \diam(B)R^{-n-1}<c\cdot h(\alpha)\le
\diam(B)R^{-n}\}.
\end{equation}
\item Define $K_0:= \{B\}$. This constitutes the 0'th layer for generalized Cantor
construction.
\item On step $n$ $(n\in \NN)$ we start with a collection $K_{n-1}$ of
balls. Define
$$L_n:= \bigcup_{b\in K_{n-1}} \SSS(b,R).
$$
Then, remove every ball from $L_n$ that intersects
$\Delta(R_\alpha, c\cdot h(\alpha))$ for at least one $R_\alpha\in C(n)$.
Denote by $K_n$ the collection balls that survive.
\item Finally, construct
$$
K_\infty= K_\infty(R):= \bigcap_{n=0}^\infty \bigcup_{B\in K_n} B.
$$
\end{enumerate}
By the construction $K_\infty$ is surely $(B,R,\vr)$-Cantor for some
parameter $\vr$. At the moment we do not have any restrictions on
the values of $\vr$, so theoretically $\vr_{n,n}$ could
equal $R$ and $K_\infty=\emptyset$. We must impose some
conditions on a pair $(\cR,h)$ in order to produce non-trivial generalized Cantor
sets. Note that $K_\infty$ can be constructed for all
(sufficiently large) values~$R$.

Assume next that every class $C(n)$ can be further split into subclasses
$C(n,m)$, $1\le m\le n$ such that for every ball $b\in K_{n-m}$ we
have
\begin{equation}\label{eq6}
\#\{D\in \SSS(b,R^m)\cap L_n\;:\; \exists R_\alpha\in C(n,m),\;
D\cap \Delta(R_\alpha,c\cdot h(\alpha))\neq\emptyset\}\ll
f(R)^{m(1-\epsilon_0)},
\end{equation}
where $0<\epsilon_0<1$ is some absolute constant. Then, one can make
Step 4 of the above algorithm more specific:
\begin{itemize}
\item[$4.1$.]Remove every ball from $L_n$ which intersects with
$\Delta(R_\alpha,c\cdot h(\alpha))$ for at least one $R_\alpha\in C(n,1)$.
By~\eqref{eq6} it will remove at most $C_1f(R)^{1-\epsilon_0}$ balls
from each set $\SSS(b,R)$, $b\in K_{n-1}$. Here $C_1$ is some
absolute positive constant.
\item[$4.m$.] (for $1<m\le n$). In general, for each $m\in \{2,\ldots,n\}$ remove every ball from $L_n$ that intersects
$\Delta(R_\alpha,c\cdot h(\alpha))$ for at least one $R_\alpha\in
C(n,m)$. By~\eqref{eq6} this process will remove at most
$C_1f(R)^{m(1-\epsilon_0)}$ balls from each set $\SSS(b,R^m)$, $b\in
K_{n-m}$.
\end{itemize}
This updated procedure ensures that $K_\infty$ is a
$(B,R,\vr)$-Cantor set with $r_{m,n}$ satisfying~\eqref{eq5} for
every $\epsilon<\epsilon_0$ and $R$ large enough.

Finally, we establish that each set $K_\infty(R)$ produced using the
bad to Cantor construction lies inside $\bad(\cR,h)$. By the
construction of each $K_n$ we have
$$\Delta(R_\alpha,c\cdot h(\alpha)) \:\: \: \cap\: \: \:  \bigcap_{n=0}^m
\bigcup_{B\in K_n} B  \: = \: \emptyset$$ for every $R_\alpha \in
\bigcup_{n=1}^m C(n)$. By letting $m$ tend to infinity we find
$K_\infty(R)\subset \bad(\cR,h)$. In other words for each $R$ large
enough there exists a $(B,R,\vr)$-Cantor subset $K_\infty(R)$ of
$\bad(\cR,h)$ with $r_{m,n}$ satisfying~$\eqref{eq5}$. This in turn
implies that $\bad(\cR,h)$ is $\epsilon_0$-Cantor winning.


To summarize, we have proved  the following theorem.

\begin{theorem}\label{th_badcant}
Let $(X,\SSS,U,f)$ be a splitting structure on $X$ and
$\bad(\cR,h)\subset X$ be a generalized bad set. If\, $\bad(\cR,h)$
adopts a bad to Cantor set construction with Condition~\eqref{eq6}
satisfied for some $\epsilon_0>0$ and some $B\in\BBB(x)$ then it is
$\epsilon_0$-Cantor-winning on $B$. In particular if $X$ satisfies
property~(S4) then
$$\dim (\bad(\cR,h)\cap
A_\infty(B)) = \dim A_\infty(B).
$$

Moreover, if the former conditions are satisfied for any ball
$B\in\BBB(X)$ and fixed $\epsilon_0>0$ then $\bad(\cR,h)$ is
$\epsilon_0$-Cantor-winning. If additionally $X$ satisfies
Property~(S5) and for any $B\in\BBB(x)$ we have $A_\infty(B) = B$ then for
any bi-Lipshitz homeomorphism $\phi\;:\; X\to X$,
$\phi(\bad(\cR,h))$ is also $\epsilon_0$-Cantor-winning.
\end{theorem}

Theorem~\ref{th_badcant} is in some sense quite cumbersome. One
needs to go through the whole bad to Cantor set construction in
order to check its conditions. In particular, one needs to construct
the sets $K_n$ and $L_n$. However, via a minor sacrifice in
generality one can improve the accessibility of
Theorem~\ref{th_badcant} and make it independent of any particular
bad to Cantor set construction. Moreover, one can ensure the
conditions are independent on the particular splitting structure. This
potentially provides the means to simultaneously establish a
Cantor-winning property of a set $\bad(\cR,h)$ for various splitting
structures $(X,\SSS,U,f)$ of the metric space $X$.

We first require some further notation. For some $R$ and $c$
satisfying \eqref{eq_rc}, assume we are given a class $C(n)$ defined
by \eqref{eq_cn} and a collection of subclasses $C(n,m)$ for $1\le
m\le n$. For any ball $b\in \BBB(X)$ let
$q_{n,m}(b)$ denote the
maximum number of balls $D\subset b$ of radius $\rad(b)R^{-m}$ such that they may intersect only on their boundaries and there exists $
R_\alpha\in C(n,m)$ satisfying $D\cap \Delta(R_\alpha,c\cdot
h(\alpha))\neq\emptyset$. Then, define
$$
q_{n,m}:= \sup\{q_{n,m}(b)\;:\; b\in\BBB(X), \rad(b) =
\rad(B)R^{m-n}\}.
$$
We may now introduce our simplification of Theorem~\ref{th_badcant}.

\begin{corollary}[C1]
Fix $B\in \BBB(X)$ and let the parameters $R$ and $c$ satisfy
\eqref{eq_rc}. Also, assume that for $n \in \NN$ we have classes
$C(n)$ defined by \eqref{eq_cn}, each associated with a collection
of subclasses $C(n,m)$ for $1\le m\le n$. If for all pairs $m,n$ and
for some $\epsilon>0$ a  splitting structure $(X,\SSS,U,f)$
satisfies $q_{n,m}\le R^{m(1-\epsilon)}$, then $\bad(\cR, h)$ is
$\epsilon$-Cantor-winning on $B$ with respect to $(X,\SSS,U,f)$. In particular, we have
$$
\dim(\bad(\cR,h)\cap A_{\infty}(B)) = \dim A_\infty(B).
$$
\end{corollary}

\proof We must simply apply the bad to Cantor set construction.
Every ball $b$ in $K_{n-m}$ has radius $\rad(b) = \rad(B)\cdot
R^{m-n}$ and all balls in $\SSS(b,R^m)\cap L_n$ are disjoint and
have radius $\rad(b)R^{-m}$. Therefore, the expression on l.h.s.
of~\eqref{eq6} does not exceed $q_{n,m}(b)$. In turn, this is at
most $q_{n,m}$. Finally, the inequality $q_{n,m}\le
R^{m(1-\epsilon)}$ assures that \eqref{eq6} is satisfied and so
Theorem~\ref{th_badcant} can be readily applied.
\endproof


As a conclusion, Theorem~\ref{th_badcant} and its corollary suggest
the following procedure to check the Cantor-winning property for a
given set $\bad(\cR,h)$.
\begin{itemize}
\item Take any large enough $R\in\NN$ and a small fixed $c>0$ in such a way that~\eqref{eq_rc} is satisfied.
\item Construct the classes $C(n)$ defined by~\eqref{eq_cn}. This constitutes a more or less straightforward task. Then, split each $C(n)$ into suitable subclasses $C(n,m)$. This often proves trickier. For `classical' sets
$\bad(\cR,h)$ it is often sufficient to take $C(n,1) := C(n)$ and
$C(n,m):=\emptyset$ for $m\ge 2$. However, for various more `modern' badly approximable sets more care is needed in the dividing process.
\item Compute an upper estimate for $q_{n,m}$; i.e., for each small
ball $b$ of radius $\rad(B)R^{-n+m}$ consider the set
$$
\bigcup_{R_\alpha\in C(n,m)} \Delta(R_\alpha,c\cdot h(\alpha)) \cap b
$$
and estimate the number of disjoint balls of smaller radius
$\rad(B)R^{-n}$ that may intersect this set.
\item If this estimate is tight enough so that $q_{n,m}\ll f(R)^{m(1-\epsilon)}$ then for the splitting
structure $(X,\SSS,U,f)$ the set $\bad(\cR,h)$ is
$\epsilon$-Cantor-winning.
\end{itemize}

We give examples of how this procedure may be implemented for
various badly approximable sets in the next section. Beforehand, we
conclude this section with a treatment of the special case that
$\cR$ consists only of points. Here,  every $\Delta(R_\alpha, c\cdot
h(\alpha))$ is simply a ball of radius $c\cdot h(\alpha)$ and the
conditions one must check to establish the Cantor-winning property
of a set become even simpler.

By the definition of class $C(n)$, for each $R_\alpha\in C(n)$ and
for each ball $D$ of radius $\rad(D) = \rad(B)R^{-n}$ one has
$$ \rad(\Delta(R_\alpha,
c\cdot h(\alpha)))\le \rad(D).
$$
Therefore, if $X$ satisfies condition (S4) then the ball
$\Delta(R_\alpha, c\cdot h(\alpha))$ can intersect at most $C(X)$
such disjoint balls $D$ of radius $\rad(B)R^{-n}$. Hence we have
$q_{n,m}(b)\le C(X)\tilde{q}_{n,m}(b)$ where
\begin{equation}\label{eq6n}
\tilde{q}_{n,m}(b):=\#\{R_\alpha\in C(n,m)\;:\; b\cap
\Delta(R_\alpha,c\cdot h(\alpha))\neq\emptyset\},
\end{equation}
and let
$$
\tilde{q}_{n,m}:= \sup\{\tilde{q}_{n,m}(b)\;:\; b\in\BBB(X), \rad(b)
= \rad(B)R^{m-n}\}.
$$
We have proved the following corollary.
\begin{corollary}[C2]
Fix $B\in \BBB(X)$ and let the parameters $R$ and $c$ satisfy
\eqref{eq_rc}. Also, assume that for $n \in \NN$ we have classes
$C(n)$ defined by \eqref{eq_cn}, each associated with a collection
of subclasses $C(n,m)$ for $1\le m\le n$.
If for all pairs $m,n$ and for some $\epsilon>0$ a splitting
structure $(X,\SSS,U,f)$ satisfies $\tilde{q}_{n,m}\le
R^{m(1-\epsilon)}$ then $\bad(\cR, h)$ is $\epsilon$-Cantor-winning
on $B$  with respect to $(X,\SSS,U,f)$.
\end{corollary}
The values $\tilde{q}_{n,m}(b)$ are usually easier to compute than
$q_{n,m}(b)$. However, in some cases $\tilde{q}_{n,m}$ may become
much larger than $q_{n,m}$ and so Corollary C2 will not be
applicable. In that case, we will have to appeal to Corollary~C1. As
an example, we will encounter this phenomenon when we consider the
standard set $\bad_N$ and the $p$-adic set $\bad_p(N)$ in the following section. The conditions in both Corollaries C1 \& C2 should be compared with the conditions required in Theorem KTV

%
\section{Applications}\label{sec_app}

\subsection{Classical badly approximable points}\label{badex}

We start this section with the model example of classical set
$\bad_N$ of  $N$-dimensional badly approximable points and describe how one can show it is Cantor-winning. As
previously mentioned the set $\bad_N$ can be written in the form of
a generalized bad set with
$$\cR =
\{\vp/q\;:\; \vp\in\ZZ^N, q\in\NN, \gcd(p_1,\ldots,p_N,q)=1\}
$$
and $h(\vp/q) = q^{-1-1/N}$. Let $B$ be the unit box $B = [0,1]^N$.
Next, choose an arbitrarily large $R$ and a real number  $c$ such
that $c^{-1} > \sqrt[N]{N!}3R^{2}$ (since $h(\vp/q)\le 1$ the
condition $\sup_\alpha(c\cdot h(\alpha))\le 1$ is satisfied). It
follows that
\begin{equation}\label{eq15}
C(n):= \{\vp/q\in\cR\;:\; cR^n\le q^{1+1/N}< cR^{n+1}\}.
\end{equation}

Now, fix a ball $b$ of $\rad(b) = R^{-n+1}$. If
$\Delta(R_{\alpha},c\cdot h(\alpha))$ intersects $b$ then we have
$$
\vd(\cent(b), R_{\alpha}) \le \rad(b) + c\cdot h(\alpha) <3/2R^{-n+1}.
$$
In other words, every element $R_\alpha\in C(n)$ for which the
neighbourhood $\Delta(R_\alpha, c\cdot h(\alpha))$ intersects $b$
must lie inside a ball of diameter $3R^{-n+1}$ centred at
$\cent(b)$.

Assume that there are at least $N+1$ points
$R_{a_1},\ldots,R_{\alpha_{N+1}}$ such that their neighbourhoods
$\Delta(R_{\alpha_i},c\cdot h(\alpha_i))$ intersect $b$. We compute
the volume of the simplex with vertices at points
$R_{\alpha_1},\ldots,R_{\alpha_{N+1}}$. On one hand this volume must
be less than $3^NR^{-N(n-1)}$ since every vertex lies inside some
box of side length $3R^{-n+1}$. On the other hand the volume is
either zero or is bounded below by $(N!\cdot q_1q_2\cdots
q_{N+1})^{-1}$ where $R_{\alpha_i} = \vp_i/q_i$. By~\eqref{eq15} we have
$$
\frac{1}{N!}(q_1q_2\cdots q_{N+1})^{-1}\ge
\frac{1}{N!}c^{-N}R^{-N(n+1)}
> 3^NR^{-N(n-1)},
$$
which is impossible. Therefore, the area of the simplex must be
zero; in other words, all points
$R_{\alpha_1},\ldots,R_{\alpha_{N+1}}$ must lie on some
$(N-1)$-dimensional affine hyperplane in $\RR^N$. If there are less
than $N+1$ points $R_{\alpha}\in C(n)$ with $\Delta(R_\alpha, c\cdot
h(\alpha))\cap b\neq\emptyset$ then we can easily find a hyperplane
containing all of them.

The upshot is that for each $b$ of radius $R^{-n+1}$ there exists a
hyperplane $\HHH_b$ which contains all the points $R_\alpha\in C(n)$
such that $\Delta(R_\alpha, c\cdot h(\alpha))\cap b\neq \emptyset$.
Thus, define $C(n,1):=C(n)$ and for $2\le m\le n$,
$C(n,m):=\emptyset$ and consider the set
$$
\EEE_b := \{E\in \BBB(\RR^k)\;:\; E\subset b,\;\rad(E) = R^{-n},\;
\exists R_\alpha\in C(n),\; E\cap
\Delta(R_\alpha,c\cdot h(\alpha))\neq\emptyset\}.
$$
It follows that every $E\in\EEE_b$ must intersect the
$c\cdot h(\alpha)$-neighbourhood of $\HHH_b$. By construction we have that
$q_{n,m}(b)$ represents the maximal number of disjoint balls in $\EEE_b$.
The definition of $C(n)$ yields that $c\cdot h(\alpha)\le R^{-n}$, and so $q_{n,1}(b) \ll R^{N-1}$. Furthermore, we have $q_{n,1}\ll R^{N-1}$. Note that for $m\ge 2$, the value of $q_{n,m}$ is surely zero.

Consider an arbitrary splitting structure $(\RR^N, \SSS, U, f)$. By Corollary~\ref{corol1}, if $d = \dim A_\infty(B)>N-1$
then $f(R) = R^d$ and
$$
q_{n,1} \ll f(R)^{\frac{N-1}{d}} = f(R)^{1 - \frac{d-N+1}{d}}.
$$
This verifies the conditions of Corollary~C1 for $\epsilon =
\frac{d-N+1}{d}$ and therefore $\bad_N$ is
$\frac{d-N+1}{d}$-Cantor-winning for $(\RR^N, \SSS,U,f)$. In
particular for the canonical splitting structure of $\RR^N$,
$\bad_N$ is $1/N$-Cantor winning. This straightforwardly implies the following proposition.
\begin{proposition}
The set $\bad_N$ has full Hausdorff dimension; i.e., $\dim\bad_N=N$.
Moreover, if for some splitting structure of $\RR^N$ one has $\dim
A_\infty(B)>N-1$ then
$$
\dim (\bad_N\cap A_\infty(B)) = \dim A_\infty(B).
$$
\end{proposition}
For a large collection of sets $A_\infty(B)$ this result  is not
new. For example for $A_\infty(B) = B$ this is simply the classical
Jarnik theorem. Many other cases are covered by the general
framework in~\cite[Theorem 8]{kristensen_thorn_velani_2006}
discussed in the previous section. However, a construction of
Cantor-winning sets for more complicated generalized bad sets
provides the answers to some open problems.

\subsection{$p$-adically badly approximable numbers}

In this subsection we demonstrate that our broad framework allows us to prove new results in spaces different to $\RR^N$. To be precise, we consider
the set $\bad_p(N)$ of $p$-adically badly approximable vectors.

\begin{theorem}\label{th_badp}
    The set $\bad_p(N)$ is $\frac{d-N+1}{d}$-Cantor-winning on any given ball $B \in \BBB(\ZZ_p^N)$ for any non-trivial splitting structure of~$\ZZ_p^N$ satisfying $d=\dim A_\infty(B)>N-1$. In particular, the set $\bad_p(N)$ is $\frac{1}{N}$-Cantor-winning with respect to the canonical splitting structure of $\ZZ^N_p$
    induced from $\QQ^N_p$.
\end{theorem}

In particular, for the canonical splitting structure of $\ZZ^N_p$
Theorem~\ref{th_badp} shows that the set $\bad_p(N)$ has maximal
Hausdorff dimension $N$, reproducing the results of~\cite{Aber}
and~\cite{kristensen_thorn_velani_2006}. However, to the best of the
authors' knowledge, no winning-type results for $\bad_p(N)$ was
previously known.

\proof

Note that the radius of any ball in $\QQ_p^N$ is an integer power of
$p$. Therefore without loss of generality we will assume that in the
proof the parameter $R$ is always an integer power of~$p$. Recall
that $\ZZ_p^N$ comes equipped with a normalized Haar measure $m$
such that the measure of each ball $b$ is $m(b) = (\rad(b))^N$.

As discussed earlier, one can readily verify that the set $\bad_p(N)$ is a generalised badly approximable
set with
$$
\cR=\{\mathbf{r}/q\in\ZZ_p^N\;:\; \mathbf{r}=(r_1, \ldots, r_N) \in
\ZZ^N, \; q\in \NN\}
$$
and $h(\mathbf{r}/q)= (\max\{|r_1|, \ldots, |r_N|,|q|\})^{-
\frac{N+1}{N}} \cdot|q|_p^{-1}$. For simplicity we provide the proof
for the particular ball $B = \ZZ^N_p$, the proof for other balls
follow the same arguments. We therefore assume from here on that
$\diam(B) = 1$ and that $q$ is always coprime with $p$ which
simplifies the formula for the height:
$$
h(\vr/q) = (\max\{|r_1|, \ldots, |r_N|,|q|\})^{- \frac{N+1}{N}}.
$$
Choose an arbitrarily large $R$, which is a power of~$p$, and a
sufficiently small $c$ to be specified later. It follows that
$$
C(n) = \{\mathbf{r}/q\in \cR\;:\; c^{-1}R^{-n-1}< h(\mathbf{r}/q) \le c^{-1}R^{-n}\}.
$$

Now, fix a ball $b$ of $\rad(b) = R^{-n+1}$. If for some $\vr/q\in
C(n)$ the neighbourhood $\Delta(\mathbf{r}/q, \:c\cdot
h(\mathbf{r}/q))$ intersects $b$, then since $c\cdot h(\mathbf{r}/q)
< R^{-n+1}$ it follows from the ultra-metric inequality that this
neighbourhood must in fact be contained in $b$. In other words,
every element $\mathbf{r}/q\in C(n)$ for which $\Delta(\mathbf{r}/q,
\: c\cdot h(\vr/q))$ intersects $b$ must lie inside $b$.

Assume as in \S\ref{badex} that there are at least $N+1$ points
$\mathbf{r}^{(1)}/q^{(1)},\ldots, \mathbf{r}^{(N+1)}/q^{(N+1)}$ not
all lying on some $(N-1)$-dimensional affine subspace of $\ZZ_p^N$
and such that their neighbourhoods $\Delta(\mathbf{r}^{(i)}/q^{(i)},
\, c\cdot h(\mathbf{r}^{(i)}/q^{(i)}))$ all intersect $b$. These
$N+1$ points therefore span a $p$-adic simplex in $\ZZ^N_p$
contained in $b$. Furthermore, by Lutz \cite{lutz_1955} the Haar
measure of this simplex is non-zero and bounded below by
\begin{equation}\label{eqn:matrix}
 c_1 \cdot \left| \det \left( \begin{array}{cccc}
1 & r_1^{(1)}/q^{(1)} & \cdots & r_{N}^{(1)}/q^{(1)}  \\
1 & r_1^{(2)}/q^{(2)} & \cdots & r_{N}^{(2)}/q^{(2)}  \\
\vdots & \vdots &  & \vdots  \\
1 & r_1^{(N+1)}/q^{(N+1)} & \cdots & r_{N}^{(N+1)}/q^{(N+1)}  \\
\end{array} \right) \right|_p,
\end{equation}
where $c_1>0$ is some absolute constant depending only upon $N$. It
is easy to check using the definitions of $C(n)$ and the height
function $h$ that the above determinant takes the form of a non-zero
rational number $M/Q$ with denominator $Q= \prod_{i=1}^{N+1}
q^{(i)}$ and numerator $M$ satisfying
$$
|M| < \left( c^{\frac{N}{N+1}}R^{\frac{(n+1)N}{N+1}} \cdot
 \right)^{N+1}\cdot \#S_{N+1}
\: \leq \: (N+1)! \cdot c^N \cdot R^{(n+1)N},
$$
where $S_{N+1}$ is the symmetric group on $N+1$ symbols. Indeed,
since we are assuming $|q^{(i)}|_p = 1$ for $i=1, \ldots, N+1$ it
follows that $|Q|_p = 1$ and so the quantity \eqref{eqn:matrix} is
bounded below by
$$
c_1 \cdot |M|_p \: \ge \: c_1 |M|^{-1} \: > \: \frac{c_1}{ (N+1)!
\cdot c^N \cdot R^{(n+1)N}}.
$$
Taking $c\le  \sqrt[N]{((N+1)! )^{-1} \cdot c_1} \cdot R^{-2}$ we
reach a contradiction since the of the ball $b$, which contains the
simplex, equals $R^{-(n-1)N}$. Moreover, it is easy to see that
criterion \eqref{eq_rc} is satisfied for this choice so long as $R$
is sufficiently large.

We deduce that all of the points $\mathbf{r}^{(1)}/q^{(1)},\ldots, \mathbf{r}^{(N+1)}/q^{(N+1)}$ must lie on some
$(N-1)$-dimensional affine hyperplane in $\ZZ_p^N$. If there are less
than $N+1$ points $\mathbf{r}/q\in C(n)$ with $\Delta(\mathbf{r}/q, \, c\cdot
h(\mathbf{r}/q))\cap b\neq\emptyset$ then we can easily find a hyperplane containing all of them. Thus, setting $C(n,1):=C(n)$ and $C(n,m):=\emptyset$ for $2\le m\le n$, it follows from geometric arguments analogous to those exhibited in \S\ref{badex} that  $q_{n,1}(b) \ll R^{N-1}$. Furthermore, we have $q_{n,1}\ll R^{N-1}$ and for $m\ge 2$ that $q_{n,m}=0$. As before, this is enough to show that for any non-trivial splitting structure on $\ZZ_p^N$ with $d=\dim A_\infty(B)>N-1$ the conditions of Corollary~C1 are satisfied with $\epsilon = \frac{d-N+1}{d}$. Therefore, the set  $\bad_p(N)$ is $\frac{d-N+1}{d}$-Cantor-winning on the ball~$B$.

\endproof

%
%
%
%

\subsection{The mixed Littlewood conjecture and the behavior of the Lagrange constant for multiples of a fixed irrational number }

Recall that the Lagrange constant $c(\alpha)$ of an irrational
number $\alpha$ is defined as the quantity
$$
c(\alpha):=\liminf_{q\to\infty}q\cdot ||q\alpha||.
$$
Obviously, $c(\alpha)>0$ if and only if $\alpha\in\bad$. On the
other hand, a  classical theorem of Dirichlet in the theory of
Diophantine approximation implies that $c(\alpha)$ cannot exceed 1.
In recent years there has been a surge of interest in investigating
the behaviour of the Lagrange constant of multiples of $\alpha$;
that is,  the behaviour of the sequence of real numbers $c(n\alpha)$ for $n\in\NN$.

By denoting $q' = qn$ one can easily observe that
$$
\liminf_{q\to\infty}qn\cdot||q\alpha|| \ge
\liminf_{q'\to\infty}q'\cdot ||q'n\alpha|| = 1/n\cdot
\liminf_{q'\to\infty} q'n\cdot ||q'n\alpha||,
$$
which in turn shows that for any positive integer $n$ and any badly
approximable $\alpha$ we always have
$$\frac{c(\alpha)}{n}\le c(n\alpha)\le nc(\alpha).
$$

In~\cite{BBEK_2014} the authors posed the following problem.
\begin{problem}\label{prob1}
Is it true that every badly approximable real number $\alpha$
satisfies
$$
\lim_{n\to\infty} c(n\alpha) = 0 \,?
$$
\end{problem}

By replacing $n$ with powers of a prime number $p$ the answer to
this problem is equivalent to the well known $p$-adic Littlewood
conjecture. It is the belief of the first author that  the answer to
Problem~\ref{prob1} is negative, although at the moment this problem
remains open. The strongest related result currently found in the
literature is due to Einsiedler, Fishman \&
Shapira~\cite{einsiedler_fishman_shapira_2011}. They answered
positively a weaker version of Problem A:
\begin{theoremEFS}
Every badly approximable real number $\alpha$ satisfies
$$
\inf_{n\ge 1} c(n\alpha) = 0.
$$
\end{theoremEFS}

Using the framework layed out in this paper we can show that there are a multitude of numbers $\alpha\in \RR$ for which the sequence
$c(n\alpha)$ either  does not tend to zero or  tends to zero as slow
as you wish.
\begin{theorem}\label{th6}
For any function $g\;:\NN\;\to \RR_{\ge 0}$ such that
$\lim_{q\to\infty}g(q) = \infty$, the set of real numbers $\alpha\in
[0,1]$ satisfying the inequality
$$
\limsup_{k\to\infty} g(k)\cdot c(k\alpha) >0
$$
is 1-Cantor-winning for any non-trivial splitting structure of
$\RR$.
\end{theorem}

{\bf Remark.} It was recently pointed out in \cite{BFS} that this result answers the dimension one case of the second part Problem~$4.4$ of Bugeaud's paper \cite{bugeaud_2015}.

\proof For any function $g$ and large parameter $R$ we will provide
the sequence $(k_i)_{i\in\NN}$ of positive integers such that
$g(k_i)\cdot c(k_i\alpha) > c$ for some positive constant $c$,
possibly dependent on $\alpha$. Then one can easily see that the set
of interest
\begin{equation}\label{eq7}
\{\alpha\in\RR\;:\; \exists c>0,\forall (i,p,q)\in
\NN\times\ZZ\times\NN,\, g(k_i)\cdot q\cdot |qk_i\alpha-p|>c\}
\end{equation}
is indeed a generalized bad set with $\cR = \{R_{i,p,q,} =
p/k_iq\;:\; (i,p,q)\in \NN\times\ZZ\times\NN\}$ and $h(i,p,q) =
(g(k_i)k_iq^2)^{-1}$. The authors do not see a possibility to apply
Theorem~KTV for this set $\bad(\cR,h)$, however we will show that Corollary~C2 is applicable.

Consider the ball $B=[0,1]$, choose an arbitrary large parameter $R$
and take $c = R^{-2}$. Then choose the values $k_i$ such that
$g(k_i)\ge R^{i-1}$ for every $i\in\NN$. We can surely do this since
$g(k) \to \infty $ as $k\to\infty$. Then let
$$
C(n) = \left\{\frac{p}{k_iq}\in\cR\;:\; R^{n-2}\le
g(k_i)k_iq^2<R^{n-1}\right\}.
$$
We split the class $C(n)$ into subclasses in the following way. Set
$$C(n,m):= \{R_{i,p,q}\in C(n)\;:\; i=m\}.$$ Then, for any two
different values $p_1/k_mq_1, p_2/k_mq_2$ from the same subclass
$C(n,m)$ we have
$$
\left|\frac{p_1}{k_mq_1}-\frac{p_2}{k_mq_2}\right| \ge
\frac{1}{k_mq_1q_2} > \frac{g(k_m)}{R^{n-1}}\ge R^{-n+m}.
$$
The final inequality automatically implies that for any ball $b$ of
radius $R^{-n+m}$ we have
$$
\#\{R_\alpha\in C(n,m)\;:\; b\cap
\Delta(R_\alpha,c\cdot h(\alpha))\neq\emptyset\}\ll 1;
$$
or, in other words,  in view of~\eqref{eq6n} we have $\tilde{q}_{n,m}\ll 1$.
Thus, for any non-trivial splitting structure $(X,\SSS,U,f)$ we have
$\tilde{q}_{n,m}\ll f(R)^{1-1}$,  the conditions of Corollary~C2 are
fulfilled and so the set $\bad(\cR,h)$ exhibited in~\eqref{eq7} is
1-Cantor-winning.
\endproof

The proof of Theorem~\ref{th6} suggests that its statement remains
valid even if we make more restrictive conditions on $\alpha$.
\begin{theorem}\label{th7}
Let $g$ be a function as in Theorem~\ref{th6}. Let
$(k_i)_{i\in\NN}$ be a sequence such that
$$
\lim_{i\to \infty} \frac{g(k_{i+1})}{g(k_i)} = \infty.
$$
Then the set of $\alpha\in[0,1]$ such that
\begin{equation}\label{eq14}
\inf_{i\in\NN} g(k_i) c(k_i\alpha) >0
\end{equation}
is 1-Cantor-winning for any non-trivial splitting structure of
$\RR$.
\end{theorem}

\proof Denote by $W$ the set of $\alpha$ satisfying
condition~\eqref{eq14} and fix an arbitrary $R$. Then, there
exists a value $i_0 = i_0(R)$ such that for every $i\ge i_0$ one has
$$
\frac{g(k_{i+1})}{g(k_i)}\ge R.
$$
Next, as in the previous proof we take $B = [0,1]$, $c=R^{-2}$ and
$k'_i := k_{i_0+i}$ ($i\in\NN$), so the condition $g(k'_i)\ge
R^{i-1}$ is satisfied. Next, we split $\cR$ into classes $C(n)$ and
then into $C(n,m)$ as in the previous proof; that is, let
$$
C(n) = \left\{\frac{p}{k_iq}\in\cR\;:\; R^{n-2}\le
g(k_i)k_iq^2<R^{n-1}\right\}
$$
and
$$
C(n,m):=\{R_{i,p,q}\in C(n)\;:\; i=m\}.
$$
Finally, by following the same arguments as in Theorem~\ref{th6} we deduce that $\tilde{q}_{n,m}\ll 1$, which in turn implies that for any
non-trivial splitting structure $(X,\SSS,U,f)$ we have
$\tilde{q}_{n,m}\ll f(R)^{1-1}$. Whence, the set
$$
W_R:=\{\alpha\in\RR\;:\; \exists c>0,\forall i\in \NN,\,
g(k'_i)\cdot c(k'_i\alpha)>c\}
$$
is in fact 1-Cantor-winning for any non-trivial splitting structure
on $\RR$. Finally notice that for any $\alpha\in W_R$,
$$
\inf_{i\in\NN}g(k_i)c(k_i\alpha) = \min_{1\le i\le i_0}
\{g(k_i)c(k_i\alpha), \inf_{j>i_0} \{g(k_j)c(k_j\alpha)\}\} =
\min_{1\le i\le i_0} \{g(k_i)c(k_i\alpha), c\}>0
$$
and so each set $W_R$ is contained in $W$. This shows that $W$, as a
supset of 1-Cantor-winning set, is itself 1-Cantor-winning.
\endproof

An important application of Theorem~\ref{th7} is that it can be
applied to certain sets related to the Mixed Littlewood Conjecture
introduced in Section~\ref{sec_introbad}.
%
For a given function $g\;:\; \NN\to \RR_{\ge 0}$ and a sequence
$\DDD = (D_n)_{n\ge 0}$ we define the set
$$
\mad_\DDD(g):=\{x\in\RR\;:\; \liminf_{q\to \infty} q\cdot g(q)\cdot
|q|_\DDD\cdot ||qx||>0\}.
$$
The Mixed Littlewood Conjecture is then precisely the statement that
$\mad_\DDD(g)$ is empty when $g\equiv 1$ for any sequence~$\DDD$.
Very recently~\cite{badziahin_velani_2011}, the following result was
proven.
\begin{theoremBV}
Let $\DDD = (2^{2^n})_{n\in\NN}$. Then, the set $\mad_\DDD(g)$ has
full Hausdorff dimension for $g(q) = \log\log q\cdot \log\log\log
q$.
\end{theoremBV}

%
With help of Theorem~\ref{th7} we show that for a suitably chosen
sequences $\DDD$ one may take even slower growing function $g(q)$
than $\log\log q\cdot \log\log\log q$. In fact, one may choose a function $g(q)$  that grows arbitrarily slowly and $\mad_\DDD(g)$ is still of full Hausdorff
dimension.

\begin{corollary}[to Theorem~\ref{th7}]
Let $g\;:\NN\;\to \RR_{\ge 0}$ be a function which monotonically
tends to infinity. Then for every sequence $\DDD =
(d_i)_{i\in\NN}$ such that
$$
\lim_{i\to\infty} \frac{g(d_{i+1})}{g(d_i)}= \infty
$$
the set $\mad_\DDD(g)$ is 1-Cantor-winning.
\end{corollary}
Unfortunately, the condition $g(q)\to \infty$ is crucial for the
proof and so this corollary does not provide any counterexample
to mixed Littlewood conjecture itself. That said, the first  author does believe that the conjecture is indeed false for sufficiently rapidly growing sequences
$\DDD$.

\proof For a given function $g$ we take the sequence
$(k_i)_{i\in\NN} = \{d_i\}_{i\in\NN}$ and consider the set $W$ as
in the proof of Theorem~\ref{th7}. It follows that the set $W$ is
1-Cantor-winning. Finally, it suffices to  check that $\mad_\DDD(g)$
contains $W$. Indeed, consider $\alpha\in W$ and an arbitrary number
$q$, and let $|q|_\DDD = k_i^{-1}$. This implies that $q=k_iq'$ and
by the definition of the pseudo-norm  it immediately follows that
$q\ge k_i$. Therefore,
$$
g(q)\cdot q\cdot |q|_\DDD\cdot ||q\alpha|| \ge g(q') \cdot q'\cdot
||q' \cdot (k_i\alpha)||.
$$
The proof is complete upon application of condition~\eqref{eq14}.
\endproof

\subsection{The $\times a$, $\times b$ problem}

In his remarkable work~\cite{furstenberg_1967}, Furstenberg showed
that if $a$ and $b$ are multiplicatively independent positive
integer numbers then for every irrational $\alpha$ the set
$$
\{a^nb^m\alpha\pmod 1\;:\; n,m\in\NN\}
$$
is dense in the unit interval. Later, Bourgain, Lindenstrauss, Michel
\& Venkatesh~\cite{blmv_2009} achieved a quantitative version of
this result, which we formulate in the following way.

\begin{theoremBLMV}
Let $\Sigma:= \{a^nb^m: n,m\in\ZZ_{\ge 0}\}$ be a multiplicative
semigroup. Then for each pair $a,b$ of multiplicatively independent
integers there exists a positive constant $c = c(a,b)$ such that the
inequality
$$
||q\alpha||<(\log\log\log q)^{-c}
$$
is satisfied for infinitely many $q\in \Sigma$.
\end{theoremBLMV}

We will show that there are numbers for which
$||q\alpha||$ can not be made too small. To be precise,  given a
function $g:\NN\to \RR_{\ge 0}$ we define

$$
\bad_{\times a,\times b}(g):=\{\alpha\in\RR\;:\; \exists c>0\,\mbox{
s.t. } \forall q\in\Sigma,\,||q\alpha||\ge c\cdot(g(q))^{-1}\}.
$$

Before stating the theorem we define a  modified logarithm
function (in order to avoid the cases when $\log q=0$). Let
$$
\log^* q:=\left\{\begin{array}{rl} 1,&\mbox{if }q<e.\\
\log q,&\mbox{otherwise}.
\end{array}\right.
$$

\begin{theorem}
For any pair $a,b$ of multiplicatively independent positive integers
the set $\bad_{\times a,\times b}(g)\cap [0,1]$ is
$\frac{\epsilon}{1+\epsilon}$\,-Cantor-winning for $g(q) = (\log^*
q)^{1+\epsilon}$, where $\epsilon$ is an arbitrary positive constant.
\end{theorem}

\noindent{\bf Remark.} By using similar methods to those used
in~\cite{badziahin_velani_2011} one can show that for
$g_1(q)=\log^*q\cdot \log^*\log q$ the set $\bad_{\times a,\times
b}(g_1)\cap [0,1]$ has in fact full Hausdorff dimension. However, this would not give us the Cantor winning property for $\bad_{\times
a,\times b}(g_1)$.

\proof

As before, we first represent the set $\bad_{\times a, \times
b}(g)$ as a generalized bad set. For this reason let
$$
\cR = \left\{\frac{p}{q}\;;\; p\in\NN, q\in\Sigma\right\}
$$
and $h(p/q) = (q\cdot g(q))^{-1}$. Consider the values
$$
r(q):=\frac{g(q)}{g(q\cdot g(q))}.
$$
for $q\in\Sigma$.
Obviously one has $r(q)<1$, but on the other hand $g(q)<q$ for all
$q>q_0(\epsilon)$. For these $q>q_0(\epsilon)$ we have
$$
\frac{g(q)}{g(q\cdot g(q))}>\frac{g(q)}{g(q^2)}=
\frac{g(q)}{2^{1+\epsilon}g(q)} = \frac1{2^{1+\epsilon}}.
$$
Therefore, $r(q)$ is bounded from below by a positive constant which
depends only on $\epsilon$. Define constants $c_1=c_1(\epsilon)$ and
$c_2=c_2(\epsilon)$ such that
$$
c_1:= \min_{q\in\NN}\{1/r(q)\};\quad c_2 = \max_{q\in\NN}\{1/r(q)\}.
$$
For sufficiently large $R$ the class $C(n)$ will take the form
$$
C(n):=\{p/q\in \cR\;:\; cR^{n}\le q\cdot g(q)<cR^{n+1}\},
$$
for some constant $c$ to be specified later. It can be readily
verified that $C(n)$ is contained within the possibly slightly
larger class
\begin{equation}\label{eq9}
C^*(n):= \left\{p/q\in \cR\;:\; \frac{c_1\cdot cR^n}{g(cR^n)}\le q<
\frac{c_2\cdot cR^{n+1}}{g(cR^{n+1})}\right\}.
\end{equation}
Now we split $C^*(n)$ into subclasses $C^*(n,s)$ in the following
way (note that these are not the subclasses $C(n,m)$ from the bad to Cantor
set construction). Let
$$
C^*(n,s):=\{p/(a^sb^t ) \in C^*(n)\;:\; t\in\ZZ_{\ge 0}\}.
$$
It is certainly the case that $s$ is bounded below by zero. On the other hand, by~\eqref{eq9} we have that
$$
s\log a\le \log\frac{c_2c\cdot R^{n+1}}{g(cR^{n+1})}.
$$
By choosing $c$ small enough we can guarantee that $s\le \frac{\log
R}{\log a}\cdot n-1$. This means that for a fixed $n$ there are at
most $\frac{\log R}{\log a}\cdot n$ various non-empty classes
$C^*(n,s)$.

Next, consider two different elements $p_1/(a^{s}b^{t_1})$ and
$p_2/(a^{s}b^{t_2})$ from $C^*(n,s)$. We have
$$
\left|\frac{p_1}{a^sb^{t_1}} - \frac{p_2}{a^sb^{t_2}}\right| \ge
\frac{1}{a^sb^{\max\{t_1,t_2\}}} \stackrel{\eqref{eq9}}>
\frac{g(cR^{n+1})}{c_2\cdot cR^{n+1}}.
$$
By taking $c$ small enough we can guarantee that the distance
between two neighbouring numbers from $C^*(n,s)$ is at least
$g(R^n)R^{-n+2\epsilon}$.

For convenience denote $k:=R^{2\epsilon}$, and let $m$ be the minimal positive integer satisfying
\begin{equation}\label{eq12}
R^m\ge Rk\cdot g(R^n).
\end{equation}
Then, it is surely the case that $R^m<R^2k\cdot g(R^n)$ and so
\begin{equation}\label{eq11}
(n\log R)^\epsilon >
\left(\frac{R^m}{R^2k}\right)^{\frac{\epsilon}{1+\epsilon}}.
\end{equation}
There must exist a natural number  $n_0(R)$ (or more exactly $n_0(R,\epsilon)$)
such that for each $n\ge n_0(R)$ the value $m$ is no
larger than $n$. Again, by choosing $c$ small enough we are able to guarantee
that $C(n) = \emptyset$ for all $n<n_0(R)$ and therefore can assume here-on that $n\ge n_0(R)$.

Consider any ball $b$ of radius $R^{-n+m}$. Then,
$$
\#\{R_\alpha\in C^*(n,s)\;:\; b\cap
\Delta(R_\alpha,c\cdot h(\alpha))\neq\emptyset\}\le
\frac{\diam(b)}{kg(R^n)R^{-n}}+2.
$$
Since $\diam(b) = R^{-n+m}$, and by the choice of $m$, the first
summand on the r.h.s. is at least $R$ and therefore for $R$ large
enough (namely $R\ge 3$) we have that the r.h.s is bounded above by
$$
\frac{2\diam(b)}{kg(R^n)R^{-n}}.
$$
Now, by collecting all classes $C^*(n,s)$ together we have
\begin{eqnarray}
\#\{R_\alpha\in C(n)\;:\; b\cap
\Delta(R_\alpha,c\cdot h(\alpha))\neq\emptyset\}\!\!\!\!\!\!\!\!\!&&\le
\frac{2R^m}{k(n\log
R)^{1+\epsilon}} \cdot \frac{n\log R}{\log a}\nonumber\\[1.5ex] \label{eq13}
&&\stackrel{\eqref{eq11}}<
\frac{2(R^2k)^{\frac{\epsilon}{1+\epsilon}}}{k\log a}\,
R^{m(1-\frac{\epsilon}{1+\epsilon})} \ll
R^{m(1-\frac{\epsilon}{1+\epsilon})}.
\end{eqnarray}
%
%

Finally,  we are ready to split $C(n)$ into subclasses to finish the
proof. Define $C(n,m)$ to be the empty set for every $n<n_0(R)$ and for every
$m\neq m_0$ for $m_0$ given by~\eqref{eq12}. Let $C(n,m_0)
= C(n)$. Then, inequality~\eqref{eq13} implies that $\tilde{q}_{n,m}\ll
R^{m(1-\epsilon/(1+\epsilon))}$ and application of Corollary C2 yields that the set
$\bad_{\times a \times b}(g)$ is
$\frac{\epsilon}{1+\epsilon}$-Cantor-winning.
\endproof

\subsection{Further examples}

In several recent papers constructions similar to generalized Cantor
sets were made inside other sets falling into the category of
generalized bad sets. With a bit of effort one can prove a
Cantor-winning property for the sets in question.

{\bf The set of points in $\bad(i,j)$ lying on vertical lines.}

Consider the pair $(i,j)$ of non-negative real numbers such that
$i+j=1$. Let $\L_x$ be a vertical line passing through the point
$(x,0)$, where $x$ satisfies the condition
\begin{equation}\label{eq16}
\liminf_{q\to \infty} q^{1/i}\cdot ||qx|| >0.
\end{equation}
To proof of Schmidt's conjecture
in~\cite{badziahin_pollington_velani_2011} the authors essentially applied a generalized Cantor set construction. To be exact, Theorem~4 and statement~(26) from~\cite{badziahin_pollington_velani_2011}
immediately imply the following.

\begin{propositionBPV}
The projection of $\bad(i,j)\cap \L_x$ onto $y$-axis is
$\epsilon$-Cantor-winning for $\epsilon = \frac1{32} (ij)^2$.
\end{propositionBPV}

Once Proposition~BPV is established one can immediately
prove a result concerning the non-empty intersection of sets
$\bad(i,j)$ for various pairs $(i,j)$. This was essentially the statement of Schmidt's conjecture.

\begin{theoremBPV}
Let $((i_\alpha,j_\alpha))_{\alpha\in S}$ be a sequence of
pairs of positive real numbers indexed by a finite or countable set $S$ such
that $i_\alpha+j_\alpha = 1$. Define
$$
i:=\inf\{i_\alpha\;:\; \alpha\in S\}\quad\mbox{ and }\quad
\epsilon:=\inf\left\{\frac{1}{32} (i_\alpha j_\alpha)^2\;:\;
\alpha\in S\right\}.
$$
Assume that $\epsilon>0$. Then for every $x\in\RR$
satisfying~\eqref{eq16}, the projection of
$$
\bigcap_{\alpha\in S} \bad(i_\alpha,j_\alpha) \cap \L_x
$$
onto $y$ axis is $\epsilon$-Cantor-winning.
\end{theoremBPV}

{\bf Sets $\bad(i_1,i_2,\ldots,i_N)$ on non-degenerate curves}

Later, in~\cite{badziahin_velani_2013}, the authors demonstrated that a result similar to Theorem~BPV holds for the sets $\bad(i,j)\cap \CCC$ for any
non-degenerate planar curve. Independently~\cite{beresnevich_2013}
Beresnevich proved more general result in higher dimensions:

Let a curve $\CCC$ be parameterized by a map
$$
\vf\;:\; I\to \RR^N;\quad \vf\in C^n(I),
$$
where $I\subset \RR$ is some interval. We assume that $\vf$ is
non-degenerate at every point on $I$ or, equivalently, that the
Wronskian of $f_1',\ldots,f_n'$ is not zero at every point $x\in
I$. Let $i_1,i_2,\ldots,i_N$ be positive real numbers such that
$i_1+\ldots+i_N=1$. Proposition~3 from~\cite{beresnevich_2013}
implies the following.
\begin{theoremB}
The set
$$
\{x\in I\;:\; \vf(x)\in \bad(i_1,\ldots,i_N)\}
$$
is $\epsilon$-Cantor-winning where
$$
\epsilon = \min\left\{(2N)^{-4}, \frac{1-(1+\min\{i_k\;:\; 1\le k\le
N\})^{-1}}{2}\right\}.
$$
\end{theoremB}

Theorem~B immediately gives a positive answer to a problem
raised by Davenport: that there are uncountably many points from
$\bad(i_1,\ldots,i_N)$ on any non-degenerate curve. In fact, the set
of such points has full Hausdorff dimension. Moreover with some
effort (see the section of \cite{beresnevich_2013} entitled `Theorem 2 implies Theorem 1') Theorem~B implies that the dimension of points from
$\bad(i_1,\ldots,i_N)$ on any non-degenerate manifold $\MMM$ is of
full Hausdorff dimension; i.e.
$$
\dim (\bad(i_1,\ldots,i_N)\cap \MMM) = \dim \MMM.
$$

%
%
%
%
%


\begin{thebibliography}{99}

\bibitem{an_2013} J. An, Badziahin-Pollington-Velani's theorem and Schmidt's game, Bull. Lond. Math. Soc. 45, No. 4, 2013, pp. 721--733.


\bibitem{Aber} A.G. Abercrombie, \textit{Badly approximable $p$-adic integers}, Proc. Indian Acad.
Sci. Math. Sci., 105(2), pp. 123--134, 1995.

\bibitem{badziahin_2013}
D. Badziahin, On multiplicatively badly approximable numbers. {\it
Mathematika}, V. 59(1), pp. 31 -- 55, 2013.

\bibitem{BHN}
D. Badziahin, S. Harrap, E. Nesharim, Topological games and Cantor-winning sets. In preparation.

\bibitem{BBEK_2014}
D. Badziahin, Y. Bugeaud, M. Einsiedler, D. Kleinbock, On the
complexity of a putative counterexample to the $p$-adic Littlewood
conjecture. Preprint, 2014.

\bibitem{badziahin_pollington_velani_2011}
D. Badziahin, A. Pollington, S. Velani, On a problem in simultaneous
Diophantine approximation: Schmidt's conjecture. {\it Annals of
Math.}, V. 174, pp. 1837 -- 1883, 2011.

\bibitem{badziahin_velani_2011}
D. Badziahin, S. Velani, Multiplicatively  badly approximable
numbers and  generalized Cantor sets. {\it Advances in Math.}, V.
228(5), pp. 2766 -- 2796, 2011.

\bibitem{badziahin_velani_2013}
D. Badziahin, S. Velani, Badly approximable points on planar curves
and a problem of Davenport. {\it Math. Ann.}, V. 359(3--4), pp. 969
-- 1023, 2014.

\bibitem{beresnevich_2013}
V. Beresnevich, Badly approximable points on manifolds. Preprint,
2013, arXiv:1304.0571.

\bibitem{blmv_2009}
J. Bourgain, E. Lindenstrauss, P. Michel, A. Venkatesh, Some
effective results for $\times a$ $\times b$. {\it Ergodic Theory
Dynam. Systems}, V. 29(6), pp. 1705–1722, 2009.


\bibitem{bfkrw_2012} R. Broderick, L. Fishman, D. Kleinbock, A. Reich, B. Weiss, The set of Badly Approximable Vectors is Strongly $C^1$ Incompressible,  Mathematical Proceedings of the Cambridge Philosophical Society, V. 153(2), pp. 319--339, 2012.

\bibitem{BFS} R. Broderick, L. Fishman, D. Simmons, Decaying and non-decaying badly approximable numbers. Preprint, 2015, arXiv:1508.03734 

\bibitem{bugeaud_2015} Y. Bugeaud,
On the multiples of a badly approximable vector, Acta Arith. 168, no. 1, pp. 71--81, 2015.


\bibitem{dani_1989} S.G.Dani, On badly approximable numbers, Schmidt games and bounded orbits of flows, Number theory and dynamical systems (York, 1987), London Math. Soc. Lecture Note Ser., V. 134, Cambridge Univ. Press, Cambridge, 1989, pp. 69--86.


\bibitem{davenport_1964}
H. Davenport, A note on Diophantin approximation II, {\it
Mathematika}, V. 11, pp. 50 -- 58, 1964.

\bibitem{einsiedler_fishman_shapira_2011}
M. Einsiedler, L. Fishman, U. Shapira, Diophantine approximations on
fractals, {\it Geom. Funct. Anal.} V. 21(1), pp. 14 -- 35, 2011.

\bibitem{fishman_simmons_urbanski_2013} L. Fishman, D. Simmons,  M. Urbanski, Diophantine Approximation, the Geometry of Limit Sets in Gromov Hyperbolic Metric Spaces, Preprint, 2013, arXiv:1301.5630v9.



\bibitem{furstenberg_1967}
H. Furstenberg, Disjointness in ergodic theory, minimal sets, and a
problem in Diophantine approximation. {\it Math. Systems Theory}, V.
1, pp. 1 -- 49, 1967.

\bibitem{jarnik_1928}
I. Jarn\'ik, Zur Metrischen Theorie der Diophantischen
Approximationen. {\it Prace Mat-Fiz.}, V. 36, pp. 91 -- 106,
1928--29.

\bibitem{kleinbock_weiss_2010}
D. Kleinbock, B. Weiss, Modified Schmidt games and Diophantine
approximation with weights. {\it Adv. Math.} V. 223(4), pp. 1276 --
1298, 2010.


\bibitem{kristensen_thorn_velani_2006}
S. Kristensen, R. Thorn, S. Velani, Diophantine approximation and
badly approximable sets. {\it Advances in Math.}, V. 203, pp. 132 --
169, 2006.

\bibitem{lutz_1955} E. Lutz, {\it Sur les approximations Diophantiennes linґeaire p?adiques}. Publications de l'Institut de Mathйmatique de l'Universitй de Strasbourg XII, Herman et Cie, Paris, (1955).

\bibitem{mcmullen_2010}
C. McMullen, Winning sets, quasiconformal maps and Diophantine
approximation. {\it Geom. Func. Anal.} V. 20(3), pp. 726 -- 740,
2010.

\bibitem{dMT} B.~de Mathan \& O.~Teuli\'{e}, \textit{Problиmes Diophantiens simultan\'{e}s}, Monatsh. Math. 143 (2004), 229-245.

\bibitem{mayeda_merrill_2013} D. Mayeda, K. Merrill, Limit points badly approximable by horoballs, Geometriae Dedicata V. 163(1), pp 127--140, 2013.

\bibitem{pollington_velani_2002}
A. Pollington, S. Velani, On simultaneously badly approxaimble
numbers. {\it J. London Math. Soc.} V. 66(2), pp. 29 -- 40, 2002.

\bibitem{schmidt_1966}
W. Schmidt, On badly approximable numbers and certain games. {\it
    Trans. Amer. Math. Soc.} V. 123, pp. 178 -- 199, 1966.

\bibitem{schmidt_1980}
W. Schmidt, Diophantine Approximation,  Lecture Notes in Mathematics 785, Springer-Verlag, 1980.

\bibitem{Weil}  S.~Weil, Schmidt games and conditions on resonant sets. Preprint, 2012, arXiv:1210.1152.


\end{thebibliography}
\end{document}